\documentclass[12pt]{article}
\usepackage{amsmath,amssymb,amsthm,amscd}
\usepackage{latexsym}
\usepackage{epsfig}

\makeatletter
\newtheorem*{rep@theorem}{\rep@title}
\newcommand{\newreptheorem}[2]{%
\newenvironment{rep#1}[1]{%
 \def\rep@title{#2 \ref{##1}}%
 \begin{rep@theorem}}%
 {\end{rep@theorem}}}
\makeatother

\newtheorem{lemma}{Lemma}[section]
\newtheorem{prop}[lemma]{Proposition}
\newtheorem{theorem}[lemma]{Theorem}
\newtheorem{thma}{Theorem}

\newreptheorem{thma}{Theorem}

\theoremstyle{definition}
\newtheorem{definition}[lemma]{Definition}
\newtheorem{remark}[lemma]{Remark}
\newtheorem{corollary}[lemma]{Corollary}
\newtheorem*{cor-nonum}{Corollary}

\numberwithin{equation}{section}

\def\deg{\text{deg} \thinspace}
\def\diag{\text{diag} \thinspace}
\def\det{\text{det} \thinspace}
\def\tr{\text{tr} \thinspace}

\newcommand{\Z}{{\mathbb Z}}
\newcommand{\C}{{\mathbb C}}

\newcommand{\PP}{{\mathbb P}}

\begin{document}

\title{\vspace{-45pt}Integrability of the Pentagram Map}

\date{June 2011\\\normalsize Revised: February 2013}

\author{\large\bf Fedor Soloviev\thanks{Department of Mathematics,
University of Toronto, Toronto, ON M5S 2E4, Canada;
e-mail: {\it soloviev\_at\_math.toronto.edu}}}

\maketitle

\begin{abstract}
The pentagram map was introduced by R. Schwartz in 1992 for convex planar polygons.
Recently, V. Ovsienko, R. Schwartz, and S. Tabachnikov proved
Liouville integrability of the pentagram map for generic
monodromies by providing a Poisson
structure and the sufficient number of integrals in involution on the space of twisted polygons.

In this paper we prove algebraic-geometric integrability for any monodromy,
i.e., for both twisted and closed polygons.
For that purpose we show that the pentagram map can
be written as a discrete zero-curvature equation with a spectral
parameter, study the corresponding spectral curve, and the dynamics on its Jacobian.
We also prove that on the
symplectic leaves Poisson brackets discovered for twisted polygons
coincide with the symplectic structure obtained from
Krichever-Phong's universal formula.
\end{abstract}

\section*{Introduction}

\begin{figure}
\begin{center}
\epsfig{file=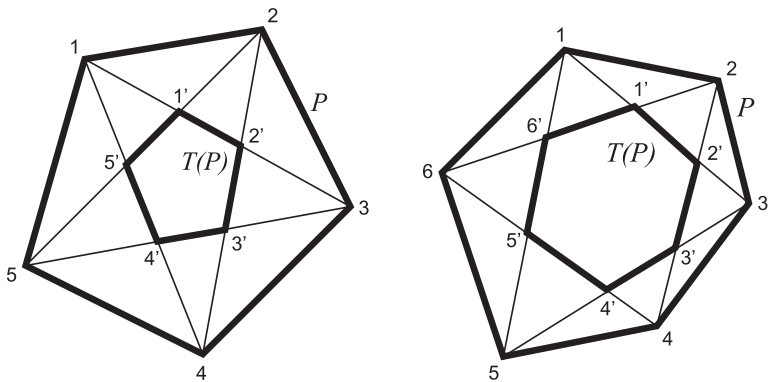,width=4in,keepaspectratio}\\
\caption{The pentagram map defined on a pentagon and a hexagon}\label{phex}
\end{center}
\end{figure}

The pentagram map was introduced by R. Schwartz in~\cite{s92} as a map defined on convex polygons understood up to
projective equivalence on the real projective plane.
Figure~\ref{phex} represents the map for a pentagon and a hexagon.

This map sends an $i$-th vertex to the intersection of 2 diagonals: $(i-1,i+1)$ and $(i,i+2)$.
The definition implies that this map is invariant under projective transformations.

Surprisingly, this simple map stands at the intersection of many branches of mathematics:
dynamical systems, integrable systems, projective geometry, and cluster algebras.
In this paper we focus on integrability of the pentagram map.

Its integrability was thoroughly studied in the paper~\cite{ost10}, where the authors considered
the pentagram map on a more general space $\mathcal{P}_n$ of the so-called twisted polygons (or $n$-gons).
A twisted polygon is a piecewise linear curve, which is not necessarily closed, but has a monodromy relating
its vertices after $n$ steps (we state its precise definition in the next section).
They proved the Arnold-Liouville integrability for the pentagram map on this space:
\begin{theorem}[\cite{ost10}]\label{ost-th}
There exists a Poisson structure invariant under the pentagram map on the space $\mathcal{P}_n$ of twisted $n$-gons.
When $n$ is even, the Poisson brackets have 4 independent Casimirs, and $n-2$ invariant functions in involution.
When $n$ is odd, there are only 2 Casimirs, and $2q$ (where $q=\lfloor n/2 \rfloor$) invariant functions in involution.
Here $\lfloor x \rfloor$ is the floor (i.e., the greatest integer) function of $x$.
\end{theorem}
The total dimension of $\mathcal{P}_n$ for all monodromies
together is $2n$, and this theorem implies the Arnold-Liouville
complete integrability on $\mathcal{P}_n$. In other words, a Zariski open
subset of $\mathcal{P}_n$ is foliated into tori, and the time
evolution is a quasiperiodic motion on these tori. The authors
of~\cite{ost10} posed an open question about integrability for
regular closed polygons. Closed polygons form a submanifold
$\mathcal{C}_n$ of codimension $8$ in $\mathcal{P}_n$, but it is
difficult to find out what happens with integrability on this
submanifold. One of the main results of the present study is a
solution of this problem in the complexified case (see Theorem~\ref{closed-poly} below).

Note that R. Schwartz conjectured that the pentagram map is a quasi-periodic motion in~\cite{s92},
introduced the integrals of motion and proved their algebraic independence in~\cite{s08}.

The central component of the algebraic-geometric integrability is a Lax representation with a spectral parameter,
which is introduced for the pentagram map in Theorem~\ref{lax-func}.
There are several advantages of this approach over the one taken in~\cite{ost10}:
\begin{itemize}
\item It works equally well both in the continuous and discrete cases. In particular,
the same algebraic-geometric methods can be used to integrate the continuous limit of
the pentagram map - the Boussinesq equation.
\item It can be used almost without changes to prove integrability for closed polygons.
\item The Lax representation provides a systematic way to obtain a Hamiltonian
structure on the space $\mathcal{P}_n$ by the universal techniques of Krichever and Phong
(more precisely, these techniques allow one to find a natural presymplectic form, which becomes
symplectic on certain submanifolds and has action-angle coordinates).
\end{itemize}

Our main results can be formulated in the following 4 theorems.
Later on we will introduce the notion of spectral data which
consists of a Riemann surface, called a spectral curve, and a
point in the Jacobian (i.e., the complex torus) of this curve.
\begin{thma}\label{spectral-th}
The space ${\mathcal P}_n$ of twisted $n$-gons (here $n \ge 4$) has a Zariski open subset
which is in a bijection with a Zariski open subset of the spectral data.
A spectral curve $\Gamma_0 \subset \C\PP^2$ is determined by complex parameters $I_j,J_j, 0 \le j \le q=\lfloor n/2 \rfloor$ as follows:
$$
R(z,k) = k^3 - k^2 \left( \sum_{j=0}^q J_j z^{j-q} \right) + k \left( \sum_{j=0}^q I_j z^{q-j} \right) z^{-n} - z^{-n} = 0.
$$
Let the normalization of $\Gamma_0$ be $\Gamma$.
For generic values of the parameters, the genus of $\Gamma$ is $g=n-2$ for even $n=2q$, and $g=n-1$ for odd $n=2q+1$.
Each torus (Jacobian $J(\Gamma)$) is invariant with respect to the pentagram map.
\end{thma}
\begin{remark}
Here and below ``generic'' means the values of the parameters from some Zariski open subset of the set of all parameters
(e.g., in this theorem ``generic parameters'' form a subspace of codimension $1$ in the space of dimension $2q+2$ as follows from Theorem~\ref{T:spec-curve}).
The bijection in the theorem is called the spectral map.
\end{remark}
Note that we consider polygons on a complex projective plane instead of a real projective plane, which
does not change any formulas for the pentagram map.

Next theorem together with the previous one establishes the algebraic-geometric integrability:
\begin{thma}\label{time-evol}
  Let $[D_{0,0}] \in J(\Gamma)$ be the point that corresponds to a generic twisted polygon at time $t=0$ after applying
  the spectral map, and $[D_{0,t}]$ be the point describing the twisted polygon at an integer time $t$.
  Then $[D_{0,t}]$ is related to $[D_{0,0}]$ by the formulas:
  \begin{itemize}
   \item when $n$ is odd,
   $$
   [D_{0,t}] = [D_{0,0} - tO_1 + tW_2] \in J(\Gamma),
   $$
   \item when $n$ is even,
   $$
   [D_{0,t}] = \left[ D_{0,0} - tO_1 + \lfloor \dfrac{1+t}{2}\rfloor W_2 + \lfloor\dfrac{t}{2}\rfloor W_3 \right],
   $$
  \end{itemize}
  provided that the corresponding spectral data remains generic up to time $t$.
  Here for odd $n$ the discrete time evolution in $J(\Gamma)$ goes along a straight line, whereas for even $n$ the evolution
  is staircase-like.

  The point $O_1 \in \Gamma$ corresponds to ($z=0$, $k$ is finite),
  and the points $W_2,W_3 \in \Gamma$ correspond to ($z=\infty,k=0$).
\end{thma}
\begin{remark}
Note that the pentagram dynamics understood as a shift on
complex tori does not prevent the corresponding orbits on the
space ${\mathcal P}_n$ from being unbounded. Indeed, these complex tori
are the Jacobians of the corresponding smooth spectral curves,
while the dynamics described above takes place for generic initial
data, i.e., for points on the Jacobians whose orbits do not
intersect special divisors (see Section 3.2). A point of a generic
orbit with an irrational shift can return arbitrarily close to
such a divisor. On the other hand, the inverse spectral map is
defined outside of these special divisors and may have poles there.
Therefore the sequence in the space ${\mathcal P}_n$ corresponding to this orbit
may escape to infinity.
\end{remark}

\begin{thma}\label{closed-poly}
For generic closed polygons the pentagram map is defined only for $n \ge 5$.
Closed polygons are singled out by the condition that $(z,k)=(1,1)$ is
  a triple point of $\Gamma$. The latter is equivalent to 5 linear relations on $I_j,J_j$:
  $$
  \sum_{j=0}^q I_j = \sum_{j=0}^q J_j = 3, \quad \sum_{j=0}^q j I_j = \sum_{j=0}^q j J_j = 3q-n, \quad
  \sum_{j=0}^q j^2 I_j = \sum_{j=0}^q j^2 J_j.
  $$
  The genus of $\Gamma$ drops to $g=n-5$ when $n$ is even, and to $g=n-4$ when $n$ is odd.
  The dimension of the Jacobian $J(\Gamma)$ drops by $3$ for closed polygons.
  Theorem~\ref{spectral-th} holds with this genus adjustment on the space $\mathcal{C}_n$, and Theorem~\ref{time-evol} holds
  verbatim for closed polygons.
\end{thma}
The relations on $I_j,J_j$ found in Theorem 4 in~\cite{ost10} are equivalent to those in Theorem~\ref{closed-poly}.
\begin{cor-nonum}
 The dimension of the phase space $\mathcal{C}_n$ in the periodic case is $2n-8$.
 In the complexified case, a Zariski open subset of $\mathcal{C}_n$ is fibred over the base of dimension $2q-3$.
 The coordinates on the base are $I_j,J_j,0 \le j \le n-1,$ subject to the constraints from Theorem~\ref{closed-poly}.
 The fibres are Zariski open subsets of Jacobians (complex tori) of dimension $2q-3$ for odd $n$, and of dimension $2q-5$ for even $n$.
 Note that the restriction of the symplectic form (which corresponds to the Poisson brackets on the symplectic leaves)
 to the space $\mathcal{C}_n$ is always degenerate, therefore the Arnold-Liouville theorem is not directly
 applicable for closed polygons. Nevertheless, the algebraic-geometric methods guarantee
 that the pentagram map exhibits quasi-periodic motion on a Jacobian.
 (Another way around this difficulty was suggested in~\cite{ost11}).
\end{cor-nonum}

Finally, the last theorem describes the relation of the
Krichever-Phong's formula with the Poisson structure of the pentagram
map. Krichever-Phong's universal formula (defined in~\cite{KP97,KP98}) applied
to the setting of the pentagram map provides a pre-symplectic 2-form on the space $\mathcal{P}_n$, see Section~\ref{kpform-sec}.
\begin{thma}\label{kp-form}
 Krichever-Phong's pre-symplectic 2-form turns out to be
 a symplectic form of rank $2g$ after the restriction to the leaves:
 $\delta I_q = \delta J_q = 0$ for odd $n$, and
 $\delta I_0 = \delta I_q = \delta J_0 = \delta J_q = 0$ for even $n$.
 These leaves coincide with the symplectic leaves of the Poisson structure found in~\cite{ost10}.
 The symplectic form is invariant under the pentagram map and coincides
 with the inverse of the Poisson structure restricted to the symplectic leaves.
 It has natural Darboux coordinates, which turn out to be action-angle coordinates for the pentagram map.
\end{thma}
We would also like to point out that there is some similarity between the pentagram map and the integrable
model~\cite{KP00} which corresponds to the $\mathcal{N}=2$ SUSY $SU(N)$ Yang-Mills theory with a hypermultiplet
in the antisymmetric representation.

\section{Definition of the pentagram map}

In this section, we give a definition of a twisted polygon, following~\cite{ost10},
introduce coordinates on the space of such polygons,
and give formulas of the map in terms of these coordinates.

\begin{definition}
A \emph{twisted $n$-gon} is a map $\phi: \Z \to \C \PP^2$, such that
none of the 3 consecutive points lie on one line (i.e., $\phi(j),\phi(j+1),\phi(j+2)$ do not lie on one line for any $j$)
and $\phi(k+n) = M \circ \phi(k)$ for any $k$. Here $M \in PSL(3,\C)$ is a projective transformation of the plane $\C \PP^2$
called the \emph{monodromy} of $\phi$.
Two twisted $n$-gons are equivalent if there is a transformation $g \in PSL(3,\C)$, such that $g \circ \phi_1 = \phi_2$.
The space of $n$-gons considered up to $PSL(3,\C)$ transformations is called $\mathcal{P}_n$.
\end{definition}
Notice that the monodromy is transformed as $M \to g M g^{-1}$ under transformations $g \in PSL(3,\C)$.
The dimension of $\mathcal{P}_n$ is $2n$, because a twisted $n$-gon depends on $2n$ variables representing
coordinates of $\phi(k), 0 \le k \le n-1$,
on a monodromy matrix $M$ ($8$ additional parameters), and the equivalence relation reduces the dimension by $8$.

There are 2 ways to introduce coordinates on the space $\mathcal{P}_n$:
If we assume that $n$ is not divisible by $3$, then there exists the unique lift of the points $\phi(k) \in \PP^2$
to the vectors $V_k \in \C^3$ provided that $\det(V_j, V_{j+1}, V_{j+2})=1$ for all $j$.
We associate a difference equation to the sequence of vectors $V_k$:
$$
V_{j+3} = a_j V_{j+2} + b_j V_{j+1} + V_j \text{ for all } j.
$$
The sequences $(a_j)$ and $(b_j)$ are $n$-periodic, i.e., $a_{j+n}=a_j, \: b_{j+n}=b_j$ for all $j$.
The monodromy is a matrix $M \in SL(3,\C)$, such that $V_{j+n} = M V_j$ for all $j$.
The variables $a_i,b_i, 0 \le i \le n-1$ are coordinates on the space $\mathcal{P}_n$ provided that $n \ne 3m$.
These coordinates are very natural, because they have a direct analogue in the
continuous KdV hierarchy.
The pentagram map is given by the formulas:
\begin{equation}\label{pent-formula}
 T^*(a_i) = a_{i+2} \prod_{l=1}^m \dfrac{1 + a_{i+3l+2} b_{i+3l+1}}{1 + a_{i-3l+2} b_{i-3l+1}}, \quad
 T^*(b_i) = b_{i-1} \prod_{l=1}^m \dfrac{1 + a_{i-3l} b_{i-3l-1}}{1 + a_{i+3l} b_{i+3l-1}}.
\end{equation}
Another set of coordinates was proposed in~\cite{ost10}. It is related to $a_i,b_i$ via the formulas:
 \begin{equation}\label{def-xy}
 x_i = \dfrac{a_{i-2}}{b_{i-2} b_{i-1}}, \qquad y_i = -\dfrac{b_{i-1}}{a_{i-2} a_{i-1}}.
 \end{equation}
 Their advantage is that they may be defined independently on $a_i,b_i$ (for any $n$) in a geometric way.
 The formulas for the pentagram map become local in the variables $x_i,y_i$,
 i.e., involving the vertex $\phi(j)$ itself and several neighboring ones:
 \begin{equation}\label{pent-formula-xy}
 T^*(x_i) = x_i \dfrac{1 - x_{i-1} y_{i-1}}{1 - x_{i+1} y_{i+1}}, \quad
 T^*(y_i) = y_{i+1} \dfrac{1 - x_{i+2} y_{i+2}}{1 - x_i y_i}.
\end{equation}
The proof of formulas~(\ref{pent-formula}) and~(\ref{pent-formula-xy}) is a direct calculation, which has been performed
  in~\cite{ost10}.\footnote{There is a typo in the formula (4.14) for $T^*(b_i)$ in~\cite{ost10}.}
Note that the pentagram map is defined only generically on $\mathcal{P}_n$ and it is not defined when a denominator
in the formulas~(\ref{pent-formula}) or~(\ref{pent-formula-xy}) vanishes.
Geometrically, it corresponds to the situation when after applying the pentagram map $3$ consecutive points of a polygon turn out to be on one line,
that is the image-polygon does not belong to the space $\mathcal{P}_n$.

 \section{A Lax representation and the geometry of the spectral curve}

 The key ingredient of the algebraic-geometric integrability is a Lax representation with a spectral
 parameter. First, we show that the pentagram map has such a representation.
 It implies the conservation of all invariant functions from Theorem~\ref{ost-th}.
 The Lax representation organizes these invariant functions in the form of the so-called spectral curve.
 We investigate some properties of the spectral curve, which are important for our purposes.

 A continuous analogue of the Lax representation is a zero-curvature equation, which is a compatibility condition for an over-determined
 system of linear differential equations (for example, see~\cite{T87} for details).
 In the discrete case, a system of differential equations becomes a system of linear difference equations
 on functions $\Psi_{i,t},i,t \ge 0,$ of an auxiliary variable $z$ (called the \emph{spectral parameter}):
 \begin{equation}\label{lin-eq}
 \begin{cases}
  L_{i,t}(z) \Psi_{i,t}(z) = \Psi_{i+1,t}(z)\\
  P_{i,t}(z) \Psi_{i,t}(z) = \Psi_{i,t+1}(z).
 \end{cases}
 \end{equation}
 The indices $i$ and $t$ are integers and represent discrete space and time variables.
 The initial polygon corresponds to $t=0$. We omit the index $t$,
 if all variables being considered in some formula correspond to the same moment of time $t$.

 It is convenient to represent several functions $\Psi_{i,t},i,t \ge 0$ and their relationship by the following diagram:
 \[\minCDarrowwidth15pt
 \begin{CD}
 \Psi_{i,t+1} @>L_{i,t+1}>> \Psi_{i+1,t+1} @>>> ... @>>> \Psi_{i+n-1,t+1} @>L_{i+n-1,t+1}>> \Psi_{i+n,t+1}\\
 @AP_{i,t}AA                @AP_{i+1,t}AA                   @.        @AP_{i+n-1,t}AA       @AP_{i+n,t}AA\\
 \Psi_{i,t}   @>L_{i,t}>>   \Psi_{i+1,t}   @>>> ... @>>> \Psi_{i+n-1,t}   @>L_{i+n-1,t}>>   \Psi_{i+n,t}
 \end{CD}
 \]
 Equations~(\ref{lin-eq}) form an over-determined system, whose compatibility condition imposes a relation
 on the functions $L_{i,t}$ and $P_{i,t}$. This relation is
 called a discrete zero-curvature equation.
 \begin{definition}
  A \emph{discrete zero-curvature equation} is the compatibility condition for system~(\ref{lin-eq}),
  which reads explicitly as:
  \begin{equation}\label{d-zcurv}
  L_{i,t+1}(z) = P_{i+1,t}(z) L_{i,t}(z) P_{i,t}^{-1}(z),
 \end{equation}
 where $L_{i,t}$ is called a \emph{Lax function}.
 \end{definition}

 \begin{theorem}\label{lax-func}
 If $n \ne 3m$, then a Lax function for the pentagram map is
 \[
 L_{i,t}(z) =
 \begin{pmatrix}
  -b_i   & 1 & 0\\
  -a_i/z & 0 & 1/z\\
    1    & 0 & 0
 \end{pmatrix} =
 \begin{pmatrix}
  0 & 0 & 1\\
  1 & 0 & b_i\\
  0 & z & a_i
 \end{pmatrix}^{-1},
 \]
 when $n=3m+1$, the corresponding function $P_{i,t}$ equals
 \[
 P_{i,t}=
 \begin{pmatrix}
  -a_i \lambda_{i-1} & 0                  & \lambda_{i-1}\\
  \lambda_{i-3}      & -a_{i+1} \lambda_i & b_{i-1} \lambda_{i-3}\\
  0                  & z \lambda_{i-2}    & 0
 \end{pmatrix},
 \text{ where } \lambda_i = \prod_{l=1}^m (1 + a_{i+3l+1} b_{i+3l}),
 \]
 and when $n=3m+2$, it equals
 \[
 P_{i,t}=
 \begin{pmatrix}
  -a_i \lambda_i \lambda_{i-1} (1+a_{i+1} b_i) & 0                  & \lambda_i \lambda_{i-1} (1+a_{i+1} b_i)\\
  \lambda_i \lambda_{i-2} (1+a_{i+1} b_i)      & -a_{i+1} \lambda_i \lambda_{i+1} (1+a_{i+2} b_{i+1}) & b_{i-1} \lambda_i \lambda_{i-2} (1+a_{i+1} b_i)\\
  0                                            & z \lambda_{i+1} \lambda_{i-1} (1+a_{i+2} b_{i+1})    & 0
 \end{pmatrix}.
\]
 For any $n$, the Lax function is
  $$
 \tilde{L}_{i,t}(z) =
 \begin{pmatrix}
   1/x_{i+2} & -1/x_{i+2} & 0\\
   1/z       & 0          & 1/z\\
   -y_{i+2}  & 0          & 0
 \end{pmatrix}=
 \begin{pmatrix}
   0       & 0 & -1/y_{i+2}\\
   -x_{i+2}& 0 & -1/y_{i+2}\\
   0       & z &  1/y_{i+2}
 \end{pmatrix}^{-1},
 $$
 with the corresponding function $\tilde{P}_{i,t}$:
 $$
 \tilde{P}_{i,t}(z) =
 \begin{pmatrix}
   1-x_{i+2}y_{i+2}                 & 0                             & 1-x_{i+2}y_{i+2}\\
   x_{i+1}y_{i+1}(1-x_{i+2}y_{i+2}) & 1-x_{i+1}y_{i+1}              & 1-x_{i+2}y_{i+2}\\
   0                                & - z y_{i+2}(1-x_{i+3}y_{i+3}) & 0
 \end{pmatrix}.
 $$
 In these formulas all variables $x_i,y_i,a_i,b_i, 0 \le i \le n-1,$ correspond to time $t$.
\end{theorem}
\begin{proof}
 The proof is to check that formulas~(\ref{pent-formula}) and~(\ref{pent-formula-xy}) are equivalent to equation~(\ref{d-zcurv})
 for our choice of the functions $L_{i,t}(z)$, $\tilde{L}_{i,t}(z)$, $P_{i,t}(z)$ and $\tilde{P}_{i,t}(z)$.
 Notice that all variables involved are $n$-periodic with respect to the index $i$.
 Here are some intermediate formulas, which appear in the proof:
\begin{itemize}
\item
 $
 \text{for } n=3m+1: \quad T^*(a_i) = a_{i+2} \dfrac{\lambda_{i+1}}{\lambda_{i-1}}, \quad T^*(b_i) = b_{i-1} \dfrac{\lambda_{i-3}}{\lambda_{i-1}}, \quad
 \dfrac{1+a_{i+1} b_i}{1+a_i b_{i-1}} \lambda_i = \lambda_{i-3},
 $
\item
 $
 \text{for } n=3m+2: \quad T^*(a_i) = a_{i+2} \dfrac{\lambda_{i+1}}{\lambda_i}, \quad T^*(b_i) = b_{i-1} \dfrac{\lambda_{i-2}}{\lambda_{i-1}}, \quad
 \dfrac{1+a_{i+3} b_{i+2}}{1+a_{i+1} b_i} \lambda_{i+2} = \lambda_{i-1}.
 $
\end{itemize}
\end{proof}
\begin{remark}\label{gauge-ab-xy}
The Lax matrices $\tilde{L}_{i,t}(z)$
and $L_{i,t}(z)$ are related by a gauge matrix $g_i = \text{diag}(1,b_i,-a_i)$:
$$
\tilde{L}_{i,t} = - \dfrac{b_{i+1}}{a_i} \left( g_{i+1}^{-1} L_{i,t} g_i \right).
$$
Note that if a proof of some theorem uses the Lax matrix $L_{i,t}$ and does not use the
``non-divisibility by 3'' condition, it will hold for $n=3m$ with $a_i,b_i$ being ``formal''
variables (i.e., not representing any polygon). However, if we switch to the variables $x_i,y_i$
using the formula above, the corresponding statement for the Lax matrix $\tilde{L}_{i,t}$
will have a real meaning, since the variables $x_i,y_i$ are defined for any $n$.
\end{remark}

 A discrete analogue of the monodromy matrix is a monodromy operator:
 \begin{definition}
 \emph{Monodromy operators} $T_{0,t},T_{1,t},...,T_{n-1,t}$ are defined as
 the following ordered products of the Lax functions:
 \begin{align*}
 &T_{0,t} = L_{n-1,t} L_{n-2,t} ... L_{0,t},\\
 &T_{1,t} = L_{0,t} L_{n-1,t} L_{n-2,t} ... L_{1,t},\\
 &T_{2,t} = L_{1,t} L_{0,t} L_{n-1,t} L_{n-2,t} ... L_{2,t},\\
 &  ...\\
 &T_{n-1,t} = L_{n-2,t} L_{n-3,t} ... L_{0,t} L_{n-1,t}.
 \end{align*}
 \end{definition}
 Similarly to the continuous case, one can define Floquet-Bloch solutions:
 \begin{definition}
  A \emph{Floquet-Bloch solution} $\psi_{i,t}$ of a difference equation $\psi_{i+1,t} = L_{i,t} \psi_{i,t}$ is an eigenvector
  of the monodromy operator: $T_{i,t} \psi_{i,t} = k \psi_{i,t}$.
 \end{definition}
 \begin{definition}\label{spec-func}
  A \emph{spectral function} of the monodromy operator $T_{i,t}(z)$ is
  $$
  R(k,z) = \hat{R}(C k,z)/C^3, \text{ where } \hat{R}(k,z) = -\det{(T_{i,t}(z) - k I)},\; C=(z^n \det{T_{i,t}(z)})^{1/3}.
  $$
  The {\it spectral curve} $\Gamma$ is the normalization of the compactification of the curve $R(k,z)=0$.
  {\it Integrals of motion} $I_j,J_j,\; 0 \le j \le q,$ are defined as the coefficients
  of the expansion
  \begin{equation}\label{spec-curve}
  R(k,z) = k^3 - k^2 \left( \sum_{j=0}^q J_j z^{j-q} \right) + k \left( \sum_{j=0}^q I_j z^{q-j} \right) z^{-n} - z^{-n}.
  \end{equation}
 \end{definition}
 \begin{remark}
 The Floquet-Bloch solutions are parameterized by the points $(k,z)$ of the spectral curve.
 Note that $C=1$ for the Lax matrix $L_{i,t}(z)$.
 However, we have $C = (-1)^n J_q/I_q \ne 1$ for the spectral function corresponding to $\tilde{L}_{i,t}(z)$.
 It is convenient to introduce the rescaling by $C$ for computational purposes.
 In particular, it makes proofs of the theorems in this section work without changes for all Lax matrices used in this paper.
 \end{remark}

 \begin{theorem}
 The coefficients $I_j,J_j,\; 0 \le j \le q,$ and the spectral curve are independent on $i$ and $t$.
 For the Lax matrix $L_{i,t}(z)$ the coefficients $I_j,J_j,$ are
  polynomials in $a_i,b_i,0 \le i \le n-1$, and they coincide with the invariants introduced in~\cite{ost10} when $n \ne 3m$.
 \end{theorem}
 \begin{proof}
 Equation~(\ref{d-zcurv}) implies that the monodromy operators satisfy the discrete-time Lax equation:
 $$
 T_{i,t+1}(z) = P_{i,t}(z) T_{i,t}(z) P_{i,t}^{-1}(z),
 $$
 i.e., monodromies $T_{i,t}$ are conjugated to each other for different $t$.
 Consequently, the function $\det{(T_{i,t}(z) - k I)}$ is independent on $t$.
 The monodromy operators $T_{i,t}(z)$ with a fixed $t$ and different $i$'s are also conjugated to each other,
 therefore $R(k,z)$ is independent on $i$.

 When $n \ne 3m$, the definition of $I_j,J_j$ in~\cite{ost10} is:
 \begin{equation}\label{old-IJ}
 \tr{(N_0 N_1 ... N_{n-1})} = \sum_{j=0}^q I_j s^{w(j)}, \qquad \tr{(N_{n-1}^{-1} ... N_1^{-1} N_0^{-1})} = \sum_{j=0}^q J_j s^{-w(j)},
 \end{equation}
 $$
 \text{where } N_j =
 \begin{pmatrix}
  0 & 0 & 1\\
  1 & 0 & b_j/s\\
  0 & 1 & a_j s
 \end{pmatrix}, \qquad
 w(j) = n+3j-3q.
 $$
 We observe that $L_j^{-1} = (g N_j g^{-1})/s$, where $g = \diag(s,s^2,1)$, if we identify $z=s^{-3}$
 (here $L_{i,t}(z) \equiv L_i(z)$). Since
 $$
  \tr{ (T_{i,t}^{-1}) } = \sum_{j=0}^q I_j z^{q-j}, \quad \tr{T_{i,t}} = \sum_{j=0}^q J_j z^{j-q}
 $$
  the invariants introduced in~\cite{ost10} coincide with our integrals of motion when $n \ne 3m$.
 \end{proof}

 We will need the explicit expressions for some of the integrals of motion for the Lax matrix $L_{i,t}(z)$ (see Proposition 5.3 in~\cite{ost10}):
 \begin{equation}\label{IJ-any-n}
 \text{for any $n \ne 3m$,}\quad I_q = \prod_{j=0}^{n-1} a_j, \quad J_q = (-1)^n \prod_{j=0}^{n-1} b_j,
 \end{equation}
 \begin{equation}\label{IJ-even-n}
 \text{for even $n\ne 3m$,}\quad I_0 = \prod_{j=0}^{q-1} b_{2j} + \prod_{j=0}^{q-1} b_{2j+1}, \quad
 J_0 = (-1)^q \prod_{j=0}^{q-1} a_{2j} + (-1)^q \prod_{j=0}^{q-1} a_{2j+1}.
 \end{equation}
 Note that if we consider $a_j,b_j$ as formal variables and use our definition of $I_j,J_j$, these formulas are valid for all $n$.

\begin{theorem}\label{T:spec-curve}
 A homogeneous polynomial $R(k,z,w)=0$ corresponding to~(\ref{spec-curve}) defines an algebraic curve $\Gamma_0$ in $\C \PP^2$.
 For generic values of the parameters $I_i,J_i$, this curve is singular only at 2 points: $(1:0:0),\; (0:1:0) \in \C \PP^2$.
 Its normalization $\Gamma$ is a Riemann surface of genus $g=2(n-q-1)$.
\end{theorem}

\begin{proof}
 A homogeneous polynomial that corresponds to equation~(\ref{spec-curve}) is
 $$
 R(k,z,w)=k^3 z^n - \sum_{j=0}^q J_j k^2 z^{n+j-q} w^{1-j+q} + \sum_{j=0}^q I_j k z^{q-j} w^{n+2+j-q} - w^{n+3}.
 $$
 The equation $R(k,z,w)=0$ defines an algebraic curve in $\C \PP^2$, which we denote by $\Gamma_0$.
 Singular points are the points where $\partial_k R = \partial_z R = \partial_w R = R = 0$.
 One can check that the only singular points with $w=0$ are the points $(1:0:0),\; (0:1:0) \in \C \PP^2$.
 Let us show that there are no singular points in the affine chart $(k:z:1)$. By Euler's theorem on
 homogenous functions, we have $k \partial_k R + z \partial_z R + w \partial_w R = (n+3)R$. Therefore,
 we have a system of 3 equations for the singular points:
 \[
 \begin{cases}
 \partial_k R = 3k^2 z^n - \sum_{j=0}^q 2 J_j k z^{n+j-q} + \sum_{j=0}^q I_j z^{q-j} = 0\\
 \partial_z R = nk^3 z^{n-1} - \sum_{j=0}^q (n+j-q)J_j k^2 z^{n+j-q-1} + \sum_{j=0}^{q-1} (q-j)I_j k z^{q-j-1} = 0\\
 R=k^3 z^n - \sum_{j=0}^q J_j k^2 z^{n+j-q} + \sum_{j=0}^q I_j k z^{q-j} - 1 = 0.
 \end{cases}
 \]
 These polynomials may have a solution in common only if $I_j,J_j,\; 0 \le j \le q,$ satisfy some non-trivial algebraic relation.
 Therefore, for generic values of the parameters $I_j,J_j$ there are no singular points in the chart $(k:z:1)$.
 For the same reason, one may assume that all branch points of $\Gamma_0$ on $z$-plane are simple, since
 the branch points of index $3$ are given by $3$ equations: $R=\partial_k R=\partial^2_k R=0$.

 According to the normalization theorem, there always exists the unique Riemann surface $\Gamma$ with a map
 $\sigma: \Gamma \to \Gamma_0$ biholomorphic away from the singular points. We will always work with
 the normalized curve $\Gamma$. The genus $g$ of $\Gamma$ is called the geometric genus of the algebraic curve $\Gamma_0$.
 To find it, we have to analyze the type of singularities of $\Gamma_0$, i.e., find the formal series solutions
 at the singular points.

\begin{lemma}\label{spec-sing}
 The singularities of the generic curve $\Gamma_0$ are as follows:
 \begin{itemize}
 \item if $n$ is even,
 the equation $R(k,z,1)=0$ has $3$ distinct formal series solutions at $z=0$:
 \begin{align*}
  O_1: \quad k_1 &= \dfrac{1}{I_q} - \dfrac{I_{q-1}}{I_q^2} z + O(z^2),\\
  O_2: \quad k_2 &= \left( \dfrac{J_0}{2}+\sqrt{J_0^2/4-I_q} \right) \dfrac{1}{z^q} + O\left(\dfrac{1}{z^{q-1}}\right),\\
  O_3: \quad k_3 &= \left( \dfrac{J_0}{2}-\sqrt{J_0^2/4-I_q} \right) \dfrac{1}{z^q} + O\left(\dfrac{1}{z^{q-1}}\right),
 \end{align*}
 and also $3$ solutions at $z=\infty$:
 \begin{align*}
  W_1: \quad k_1 &= J_q + \dfrac{J_{q-1}}{z} + O\left(\dfrac{1}{z^2}\right),\\
  W_2: \quad k_2 &= \left( \dfrac{I_0+\sqrt{I_0^2-4J_q}}{2J_q} \right) \dfrac{1}{z^q} + O\left(\dfrac{1}{z^{q+1}}\right),\\
  W_3: \quad k_3 &= \left( \dfrac{I_0-\sqrt{I_0^2-4J_q}}{2J_q} \right) \dfrac{1}{z^q} + O\left(\dfrac{1}{z^{q+1}}\right).
 \end{align*}

 \item if $n$ is odd,
 the equation $R(k,z,1)=0$ has $3$ distinct Puiseux series solutions at $z=0$:
 \begin{align*}
  O_1: \quad k_1 &= \dfrac{1}{I_q} - \dfrac{I_{q-1}}{I_q^2} z + O(z^2),\\
  O_2: \quad k_2 &= \dfrac{\sqrt{-I_q}}{z^{n/2}} + \dfrac{J_0}{2 z^{(n-1)/2}} + O\left(\dfrac{1}{z^{(n-2)/2}}\right),\\
             k_3 &= -\dfrac{\sqrt{-I_q}}{z^{n/2}} + \dfrac{J_0}{2 z^{(n-1)/2}} + O\left(\dfrac{1}{z^{(n-2)/2}}\right),
 \end{align*}
 and $3$ solutions at $z=\infty$:
 \begin{align*}
  W_1: \quad k_1 &= J_q + \dfrac{J_{q-1}}{z} + O\left(\dfrac{1}{z^2}\right),\\
  W_2: \quad k_2 &= \dfrac{1}{\sqrt{-J_q}} \dfrac{1}{z^{n/2}} + \dfrac{I_0}{2 J_q} \dfrac{1}{z^{(n+1)/2}} + O\left(\dfrac{1}{z^{(n+2)/2}}\right),\\
             k_3 &= -\dfrac{1}{\sqrt{-J_q}} \dfrac{1}{z^{n/2}} + \dfrac{I_0}{2 J_q} \dfrac{1}{z^{(n+1)/2}} + O\left(\dfrac{1}{z^{(n+2)/2}}\right).
 \end{align*}

 \end{itemize}
 If $\sigma: \Gamma \to \Gamma_0$ is a normalization of $\Gamma_0$, the singularities of $\Gamma_0$ correspond
 to the following points on $\Gamma$:
 \begin{itemize}
 \item for even $n, \qquad \sigma^{-1}(1:0:0) = \{O_2,O_3\}, \qquad \sigma^{-1}(0:1:0) = \{ W_1, W_2, W_3 \}$.
 \item for odd $n, \qquad \sigma^{-1}(1:0:0) = O_2, \qquad \sigma^{-1}(0:1:0) = \{ W_1, W_2 \}$,
 \end{itemize}
 The point $O_1 \in \Gamma$ is non-singular.
\end{lemma}
\begin{proof}
One finds the coefficients of the series recursively by substituting the series into the equation~(\ref{spec-curve}), which determines the spectral curve.
\end{proof}

 Now we can complete the proof of Theorem~\ref{T:spec-curve}. First, we
 find the number of branch points of $\Gamma$, and then
 we use the Riemann-Hurwitz formula to find the genus of $\Gamma$.

 The number of branch points of $\Gamma$ on $z$-plane equals the number of zeroes of the function:
 $$
 \partial_k R(k,z) =
 3 k^2 - 2k \left( \dfrac{J_0}{z^q} + \dfrac{J_1}{z^{q-1}} + ... + \dfrac{J_{q-1}}{z} + J_q \right) +
 \left( \dfrac{I_0}{z^{n-q}} + \dfrac{I_1}{z^{n-q+1}} + ... + \dfrac{I_{q-1}}{z^{n-1}} + \dfrac{I_q}{z^n} \right)
 $$
 with an exception of the singular points.
 The function $\partial_k R(z,k)$ is meromorphic on $\Gamma$, therefore the number of its zeroes equals
 the number of its poles.
 For any $n$, $\partial_k R$ has poles of total order $3n$ at $z=0$,
 and $\partial_k R$ has zeroes of total order $n$ at $z=\infty$.
 For even $n$ the Riemann-Hurwitz formula implies that $2-2g=6-(3n-n)$, thus the genus of $\Gamma$ is $g=n-2$.
 For odd $n$ we have $2-2g=6-(3n-n+2)$, and $g=n-1$. The difference between odd and even values of $n$ occurs
 because $O_2,W_2$ are branch points for odd $n$.
\end{proof}

\section{Direct and inverse spectral transforms}

In this section we prove Theorems~\ref{spectral-th} and~\ref{time-evol}.
\begin{definition}
{\rm
Let $J(\Gamma)$ be the Jacobian of the generic spectral curve $\Gamma$, and $[D]$ be a point in the Jacobian.
The pair consisting of the spectral curve $\Gamma$ (with marked points $O_i$ and $W_i$) and  a point $[D] \in J(\Gamma)$ is called the {\it spectral data.}
}
\end{definition}
\noindent Theorem~\ref{spectral-th} may be stated as follows:
\begin{repthma}{spectral-th}
The space ${\mathcal P}_n$ of twisted $n$-gons (here $n \ge 4$) has a Zariski open subset
which is in a bijection with a Zariski open subset of the spectral data.
This bijection is called the spectral map.
The Jacobians (complex tori) are invariant with respect to the pentagram map.
\end{repthma}
\begin{remark}\label{dim-count}
 $\Gamma$ is determined by $2q+2$ parameters: $I_j,J_j,0 \le j \le q$,
 and $J(\Gamma)$ has the dimension $g$, therefore
 the dimensions of the domain and the range of the spectral map match.
 The existence of this bijection implies functional independence of the parameters $I_j,J_j, 0\le j \le q$
 and coordinates in $J(\Gamma)$, as well as the fact that
 generically a divisor obtained after applying the spectral map is non-special.
 The independence of $I_j,J_j, 0 \le j \le q$ was proved in~\cite{ost10} by a different method.
\end{remark}

 The proof of Theorem~\ref{spectral-th} consists of two parts: the construction of the
 direct spectral transform $S$ and the construction of the inverse spectral transform $S_{inv}$.
 $S$ and $S_{inv}$ are inverse to each other on a Zariski open subset
 (however, the domain of the map $S$ (or $S_{inv}$) may be different from the range of $S_{inv}$(or $S$, respectively)).


 \subsection{Direct spectral transform.}\label{direct-st}

 Given a point in the space ${\mathcal P}_n$, we construct the spectral curve $\Gamma$ and the Floquet-Bloch
 solution $\psi_{0,0}$.  In our definition of the spectral data, $\Gamma$ has to be generic.
 Therefore, the domain of $S$ is a Zariski open subset ${\mathcal P}_n^0 \subset {\mathcal P}_n$
 that consists of those points, for which $\Gamma$ is generic.
 In what follows we always assume that $\Gamma$ is generic.
 The vector function $\psi_{0,0}$ is defined up to a multiplication by a scalar function.
 To get rid of this ambiguity, we \emph{normalize} $\psi_{0,0}$ by dividing it by the sum of its components.
 As a result, the vector function $\psi_{0,0}$ satisfies the identity: $\sum_{i=1}^3 \psi_{0,0,i} \equiv 1$
 (here the index $i$ denotes the $i$-th component of the vector $\psi_{0,0}$). Additionally, it acquires poles
 on the curve $\Gamma$. We denote the pole divisor of $\psi_{0,0}$ by $D_{0,0}$.
 The Abel map assigns a point in the Jacobian $J(\Gamma)$ of the curve $\Gamma$ to each divisor on $\Gamma$.
 We denote by $[D_{0,0}]$ the corresponding point in $J(\Gamma)$. It constitutes the second part of the spectral data
 and is used to define the map $S: {\mathcal P}_n^0 \to (\Gamma, [D])$.

 \begin{remark}
 A priori, the subset ${\mathcal P}_n^0$ may be empty. Its non-emptiness follows from the existence of the map $S_{inv}$.
 This argument is standard in the theory of algebraic-geometric integration.
 \end{remark}

 Notice that once we define the vector function $\psi_{0,0}$, all other vectors $\psi_{i,t}$ with $i,t \ge 0$
 are uniquely determined using equations~(\ref{lin-eq}) in the vector form:
\begin{equation*}
 \begin{cases}
  L_{i,t}(z) \psi_{i,t}(z) = \psi_{i+1,t}(z)\\
  P_{i,t}(z) \psi_{i,t}(z) = \psi_{i,t+1}(z).
 \end{cases}
 \end{equation*}
 The vectors $\psi_{i,t}$ with $i,t \ne 0$ are not normalized.
 In Theorem~\ref{time-evol} below we need to normalize each vector $\psi_{i,t}$,
 and we denote the normalized vectors by $\bar{\psi}_{i,t}$.
 The vectors $\psi_{0,0}$ and $\bar{\psi}_{0,0}$ are identical in this notation.
 The following proposition establishes the number of poles of the \emph{normalized}
 Floquet-Bloch solution for any values of $i,t$.
 \begin{prop}\label{deg-D}
 If the spectral curve $\Gamma$ corresponding to the Lax functions is generic,
 a Floquet-Bloch solution $\bar{\psi}_{i,t}$ is a meromorphic vector function on $\Gamma$. It is uniquely
  defined by the requirement $\sum_{j=1}^3 \bar{\psi}_{i,t,j} \equiv 1$. Generically, its pole divisor $D_{i,t}$ has degree $g+2$.
 \end{prop}
 \begin{proof}
 First of all, we show that $\bar{\psi}_{i,t}$ is a meromorphic function.
 By definition, it is a solution to the linear equation: $(T_{i,t}(z) - k)u = 0$.
 By Cramer's rule, the components of the vector $u$ are rational functions in the entries of the matrix
 $T_{i,t}(z) - k$ and, consequently, they are rational functions in $k$ and $z$.
 The normalized solution ($u$ divided by the sum of its components $u_1+u_2+u_3$)
 is also a rational function in $k$ and $z$, i.e., a meromorphic function on $\Gamma$.

 Secondly, we find the behavior of $\bar{\psi}_{i,t}$ at the branch points.
 Let the expansion of $k(z)$ at the branch point $(k_0,z_0) \in \Gamma$ be
 $k(z) = k_0 \pm k_1 \sqrt{z-z_0} + O(z-z_0)$. If we assume that $k_1=0$, then the equation $R(k,z)=0$ implies
 that $\partial_z R(k_0,z_0)=0$, i.e., the point $(k_0,z_0) \in \Gamma$ is singular.
 Since it is not possible by Theorem~\ref{T:spec-curve}, we have that $k_1 \ne 0$.
 One can check that the corresponding expansion of $\bar{\psi}_{i,t}$ at the branch point is
 $\bar{\psi}_{i,t} = v \pm w \sqrt{z-z_0} + O(z-z_0)$, where the vectors $v$ and $w$ are determined as follows:
 $$
 T_{i,t}(z_0) v = k_0 v, \qquad (T_{i,t}(z_0) - k_0)w = k_1 v, \qquad \sum_{i=1}^3 v_i = 1, \qquad \sum_{i=1}^3 w_i = 0.
 $$
 The latter equations determine $v,w$ uniquely, and they imply that $k_0$ corresponds to a Jordan block
 of the matrix $T_{i,t}(z_0)$.

 Thirdly, we find the number of the poles of $\bar{\psi}_{i,t}$.
 If $u_1+u_2+u_3=0$, then the function $\bar{\psi}_{i,t}$ may develop a pole. For generic values
 of the parameters $a_i,b_i$, we may assume that these poles are distinct from the branch points of $\Gamma$.
 Let $k_i,1\le i \le 3$ be the solutions of equation~(\ref{spec-curve}) for a fixed value of $z$.
 Then $Q_i=(k_i,z),1\le i \le 3,$ correspond to $3$ points on $\Gamma$, and we can form a matrix
 $\bar{\Psi}_{i,t}(z) = \{ \bar{\psi}_{i,t}(Q_1), \bar{\psi}_{i,t}(Q_2), \bar{\psi}_{i,t}(Q_3) \}$.
 Obviously, this matrix depends on the ordering of the roots $k_1,k_2,k_3$.
 However, an auxiliary function $F(z) = \det^2{\bar{\Psi}_{i,t}(z)}$ is independent on that ordering.
 Consequently, $F(z)$ is a well-defined meromorphic function on $\Gamma$. Generically, it is not singular
 at the points $z=0$ and $z=\infty$, which follows from Proposition~\ref{psi-divisor} below.
 One can check using the above series expansion of $\bar{\psi}_{i,t}$ that $F(z)$ has zeroes precisely
 at the branch points of $\Gamma$, and that these zeroes are simple.
 In Theorem~\ref{T:spec-curve} we found that the number of the branch points of $\Gamma$ is $\nu = 2g+4$.
 The pole divisor of $F(z)$ equals $2 \pi(D_{i,t})$.
 Consequently, we have $deg \thinspace D_{i,t} = \nu/2 = g+2$.
 \end{proof}

 In the following proposition we drop the index $t$, since all variables correspond to the same moment of time.
 \begin{prop}\label{psi-divisor}
  Generically, the divisors of the functions $\psi_{i,j}, 0 \le i \le n-1, 1 \le j \le 3$ satisfy the following inequalities:
  \begin{itemize}
  \item when $n$ is odd,
  $$
  (\psi_{i,1}) \ge -D + (1-i)O_2+(1+i)W_2, \qquad (\psi_{i,2}) \ge -D - iO_2+W_1+(1+i)W_2,
  $$
  $$
  (\psi_{i,3}) \ge -D + (2-i)O_2+iW_2;
  $$
  \item when $n$ is even,
  $$
  (\psi_{2k,1}) \ge -D -k O_2+(1-k)O_3+k W_2+(1+k)W_3,
  $$
  $$
  (\psi_{2k,2}) \ge -D -k O_2-k O_3+W_1+k W_2+(1+k)W_3,
  $$
  $$
  (\psi_{2k,3}) \ge -D +(1-k) O_2+(1-k) O_3+k W_2+k W_3,
  $$
  $$
   (\psi_{2k+1,1}) \ge -D -k O_2-k O_3+(1+k) W_2+(1+k)W_3,
  $$
  $$
  (\psi_{2k+1,2}) \ge -D -(1+k) O_2-k O_3+W_1+(1+k) W_2+(1+k)W_3,
  $$
  $$
  (\psi_{2k+1,3}) \ge -D -k O_2+(1-k) O_3+k W_2+(1+k)W_3.
  $$
  \end{itemize}
 \end{prop}

 \begin{proof}
 In this proof, we use the argument from Remark~\ref{gauge-ab-xy}.
 First note that our computation is formally valid for any $n$. Secondly, the components of
 the non-normalized vectors $\psi_i$ are proportional for the Lax matrices related by diagonal gauge matrices,
 therefore this proposition holds for such Lax matrices as well. In particular, Remark~\ref{gauge-ab-xy} implies
 that the vector $\tilde{\psi}_i$ for the Lax matrix $\tilde{L}_i$
 has the same divisor structure as $\psi_i$ for any $n$ and $i,\: 0 \le i \le n-1$.

 First we prove the inequalities of the theorem for the components of the vector function $\psi_0$.
 Then we use a permutation argument and Lemmas~\ref{omega-O1}-\ref{omega-W23} to complete the proof for the vector functions $\psi_i$ with $i>0$.
 The situation is different for even and odd $n$.

 When $n$ is even, using Lemmas~\ref{spec-sing} and~\ref{T0-asymp}, the definition of the Floquet-Bloch solution, and the normalization condition $\psi_{0,1}+\psi_{0,2}+\psi_{0,3} \equiv 1$,
 one can check that $\psi_0$ is holomorphic at the points $O_1,O_2,O_3$ and that
 $$
  \psi_0 = \begin{pmatrix} O(z)\\O(1)\\O(z) \end{pmatrix} \text{ at } O_3, \qquad \psi_0 = \begin{pmatrix} O(1)\\O(1)\\O(z) \end{pmatrix} \text{ at } O_2,
   \qquad a_0 = \lim_{Q \to O_1} \dfrac{\psi_{0,3}(Q)}{\psi_{0,1}(Q)}.
 $$
 Similarly, the expansion of $T^{-1}_{0,t}(z)$ at $z=\infty$ along with the identity $T_0^{-1} \psi_0 = k^{-1} \psi_0$, implies that $\psi(Q)$ is holomorphic at the points $W_1,W_2,W_3$ and that
 $$
 \psi_0 = \begin{pmatrix} O(1/z)\\O(1/z)\\O(1) \end{pmatrix} \text{at } W_3,\; \psi_0 = \begin{pmatrix} O(1)\\O(1/z)\\O(1) \end{pmatrix} \text{at } W_1, \;
 b_0 = \lim_{Q \to W_2} \dfrac{\psi_{0,2}(Q)}{\psi_{0,1}(Q)}, \;
 b_{n-1} = -\lim_{Q \to W_1} \dfrac{\psi_{0,1}(Q)}{\psi_{0,3}(Q)}.
 $$
 We perform a similar analysis for odd $n$:
 \begin{equation}\label{oddpsi}
 \psi_0 = \begin{pmatrix} O(\sqrt{z})\\O(1)\\O(z) \end{pmatrix} \text{ at } O_2, \quad
 \psi_0 = \begin{pmatrix} O(1)\\O(1/z)\\O(1) \end{pmatrix} \text{ at } W_1 \quad
 \psi_0 = \begin{pmatrix} O(1/\sqrt{z})\\O(1/\sqrt{z})\\O(1) \end{pmatrix} \text{ at } W_2,
 \end{equation}
 \begin{equation}\label{findab0}
 a_0 = \lim_{Q \to O_1} \dfrac{\psi_{0,3}(Q)}{\psi_{0,1}(Q)}, \qquad b_0 = \lim_{Q \to W_2} \dfrac{\psi_{0,2}(Q)}{\psi_{0,1}(Q)}, \qquad
 b_{n-1} = -\lim_{Q \to W_1} \dfrac{\psi_{0,1}(Q)}{\psi_{0,3}(Q)}.
 \end{equation}
 Notice that a cyclic permutation of indices $(n-1,n-2,...,1,0)$ changes $T_i \to T_{i+1}$ and
 $\bar{\psi}_i \to \bar{\psi}_{i+1}$. For even $n$, it also permutes
 $\bar{\psi}_i(O_2) \leftrightarrow \bar{\psi}_i(O_3)$ and
 $\bar{\psi}_i(W_2) \leftrightarrow \bar{\psi}_i(W_3)$.

 The latter happens for the following reason:
 The asymptotic expansions of $k$ at the points $O_2,O_3$ given by Lemma~\ref{spec-sing}
 contain expressions $J_0/2 \pm \sqrt{J_0^2/4-I_q}$,
 which are equal to $(-1)^q \prod_{j=0}^{q-1} a_{2j}$ and $(-1)^q \prod_{j=0}^{q-1} a_{2j+1}$, i.e.,
 a cyclic permutation of indices swaps the eigenvalues and the corresponding eigenvectors.
 Likewise, the expressions $\left(I_0 \pm \sqrt{I_0^2-4J_q}\right)/\left( 2J_q \right)$ at the points $W_2,W_3$
 are equal to $\left( \prod_{j=0}^{q-1} b_{2j} \right)^{-1}$ and $\left( \prod_{j=0}^{q-1} b_{2j+1} \right)^{-1}$
 and are also swapped.
 This observation allows us to produce formulas for the components of the vectors $\bar{\psi}_i$ from the formulas for $\bar{\psi}_0 \equiv \psi_0$.
 For example, for even $n$ we obtain:
 $$
 \bar{\psi}_{2j} = \begin{pmatrix} O(z)\\O(1)\\O(z) \end{pmatrix} \text{at } O_3, \: \bar{\psi}_{2j} = \begin{pmatrix} O(1)\\O(1)\\O(z) \end{pmatrix} \text{at } O_2, \:
 \bar{\psi}_{2j+1} = \begin{pmatrix} O(z)\\O(1)\\O(z) \end{pmatrix} \text{at } O_2, \: \bar{\psi}_{2j+1} = \begin{pmatrix} O(1)\\O(1)\\O(z) \end{pmatrix} \text{at } O_3.
 $$
 Now we can use Lemmas~\ref{omega-O1}-\ref{omega-W23} to complete the proof of the proposition.
 Consider, for example, the vector $\psi_i$ at the point $O_2$ for even $n$.
 The proof of Lemma~\ref{omega-O23} implies that
 $$
 \bar{\psi}_{2j} = \begin{pmatrix} 1 \\O(1)\\O(z) \end{pmatrix} \text{at } O_2, \; \text{ and that }
 \psi_{2j} = \begin{pmatrix} O(z^{-j}) \\O(z^{-j})\\O(z^{1-j}) \end{pmatrix} \text{at } O_2,
 $$
 $$
 \text{since }\psi_{2j}=\bar{\psi}_{2j}/f_{2j}(z), \text{ where } f_{2j}(z)= \dfrac{(-1)^j}{\prod_{k=0}^{j-1} a_{2k}} z^j + O(z^{j+1}) \text { at } O_2.
 $$
 Lemma~\ref{omega-O23} uses a different normalization of eigenvectors.
 However, we normalize only the vector $\psi_i$ with $i=0$ and do not
 normalize the vectors with $i>0$. This means that generically we still have
 $\psi_{2j}=(O(z^{-j}),O(z^{-j}),O(z^{1-j}))^T$ at the point $O_2$, which agrees with the multiplicity in the statement of the proposition.
 Note that Lemma~\ref{omega-O23} does not provide formulas for $f_{2j+1}(z)$ and we have to analyze the vectors $\psi_i$ with odd $i$ separately.
 The vector equation $\psi_{2j+1} = L_{2j} \psi_{2j}$ is equivalent to
 $$
 \psi_{2j+1,1} = \psi_{2j,2} - b_{2j} \psi_{2j,1},\qquad \psi_{2j+1,2} = (\psi_{2j,3} - a_{2j} \psi_{2j,1})/z,\qquad \psi_{2j+1,3} = \psi_{2j,1}.
 $$
 The latter equations imply that generically we have $\psi_{2j+1}=(O(z^{-j}),O(z^{-j-1}),O(z^{-j}))^T$ at the point $O_2$.
 Other cases (the points $O_1,O_2,O_3,W_1,W_2,W_3$ for both even and odd $n$) are treated similarly.
 The following formulas are used in the proof (they follow from the formulas above by using a permutation of indices, and they hold both for even and odd $n$ with the understanding that $W_2=W_3$ for odd $n$):
 \begin{equation*}
 a_i = \underset{Q \to O_1}{\lim} \dfrac{\bar{\psi}_{i,3}(Q)}{\bar{\psi}_{i,1}(Q)}, \qquad
 b_{2k} = \underset{Q \to W_2}{\lim} \dfrac{\bar{\psi}_{2k,2}(Q)}{\bar{\psi}_{2k,1}(Q)}, \qquad
 b_{2k+1} = \underset{Q \to W_3}{\lim} \dfrac{\bar{\psi}_{2k+1,2}(Q)}{\bar{\psi}_{2k+1,1}(Q)}.
 \end{equation*}
 Since $\bar{\psi}_i = f_i \psi_i$, we also have:
 $$
 a_i = \underset{Q \to O_1}{\lim} \dfrac{\psi_{i,3}(Q)}{\psi_{i,1}(Q)}, \qquad
 b_{2k} = \underset{Q \to W_2}{\lim} \dfrac{\psi_{2k,2}(Q)}{\psi_{2k,1}(Q)}, \qquad
 b_{2k+1} = \underset{Q \to W_3}{\lim} \dfrac{\psi_{2k+1,2}(Q)}{\psi_{2k+1,1}(Q)}.
 $$
 \end{proof}
 \begin{remark}
 Note that gauge transformations by non-degenerate diagonal matrices do not change the structure of the divisors
 of the non-normalized vectors $\psi_i$ given by Proposition~\ref{psi-divisor}.
 A normalization of $\psi_0$ changes the divisor $D_{0,0}$ to an equivalent one.
 Since we consider only Lax matrices related by such gauge transformations,
 Proposition~\ref{psi-divisor} holds for all of them.
 \end{remark}

 \subsection{Inverse spectral transform.}

 The construction of the map $S_{inv}$ is completely independent of the construction of $S$.
 It consists of $3$ parts (which we describe in detail below):
 \begin{itemize}
 \item We use the analytic properties of the Floquet-Bloch solution $\psi_i$ established in Proposition~\ref{psi-divisor}
 as a \emph{motivation} for the construction of $S_{inv}$.
 We assume that the domain of $S_{inv}$ consists of the generic spectral data.
 Here ``generic spectral data'' means spectral functions which may be singular only at the points $O_i,W_i$
 and divisors $[D] \equiv [D_{0,0}] \in J(\Gamma)$, such that all divisors in
 Proposition~\ref{psi-divisor} with $0 \le i \le n-1$ are non-special.
 This assumption allows us to use the Riemann-Roch theorem to reconstruct the components $\psi_{i,j},1 \le j \le 3,$
 of the vector $\psi_i$ up to a multiplication by constants.

 Since the number of the divisors in Proposition~\ref{psi-divisor} is finite, generic spectral data is determined
 by a finite number of algebraic relations in the space of spectral data, and hence it is a Zariski open subset.

 We drop the index $t$ below, because all variables correspond to the same moment of time.

 \item Given the generic spectral data and any non-zero complex number $C$, Proposition~\ref{ag-int} allows us
 to reconstruct Lax matrices $L_j',0 \le j \le n-1$:
  \[
  L'_j(z) =
  \begin{pmatrix}
   0 & 0 & c'_j\\
   d'_j & 0 & b'_j\\
   0 & e'_j z & a'_j
  \end{pmatrix}^{-1},
  0 \le j \le n-1, \qquad
  \prod_{i=0}^{n-1} c'_i d'_i e'_i = C^{-3}.
  \]

 \item Proposition~\ref{general-lax} allows us to perform the reduction
 from $L_j'$ to either $L_j$ or $\tilde{L}_j$, which completes the construction of $S_{inv}$
 (in the case of $L_j$ we set $C=1$ and $n$ cannot be a multiple of $3$;
  in the case of $\tilde{L}_j$ any $n$ is possible and we set $C = (-1)^n J_q/I_q$).
 \end{itemize}

 It will be evident from the construction that $S \circ S_{inv} = Id$ and $S_{inv} \circ S = Id$ when both maps $S$ and $S_{inv}$ are defined.
 Since they are defined on Zariski open subsets, their composition is also defined on a Zariski open subset.
 This concludes the proof of Theorem~\ref{spectral-th}.
 \begin{prop}\label{ag-int}
  Given the generic spectral data and any number $C \in \C \backslash \{ 0 \}$, one can recover a sequence of $n$ matrices:
  \[
  L'_i(z) =
  \begin{pmatrix}
   0 & 0 & c'_i\\
   d'_i & 0 & b'_i\\
   0 & e'_i z & a'_i
  \end{pmatrix}^{-1},
  0 \le i \le n-1, \qquad
  \prod_{i=0}^{n-1} c'_i d'_i e'_i = C^{-3}.
  \]
  This sequence is unique up to gauge transformations: $L'_i \to g_{i+1} L'_i g_i^{-1}$,
  where $g_i, 0 \le i \le n-1,$ are non-degenerate diagonal matrices ($g_n=g_0$).
 \end{prop}
 \begin{proof}
  The procedure to reconstruct the matrices $L'_i, 0 \le i \le n-1,$ consists of 3 steps:
  \begin{enumerate}
  \item We pick an arbitrary divisor $D$ of degree $g+2$ in the equivalence class $[D] \in J(\Gamma)$.
  \item We observe that the degree of all divisors in Proposition~\ref{psi-divisor} is $-g$.
  According to the Riemann-Roch theorem, it means that each function $\psi_{i,j}$ is determined up to
  a multiplication by a constant. We pick arbitrary non-zero constants, and thus obtain a sequence of vectors $\psi_i, 0 \le i \le n-1$.
  We define $\psi_n = C k \psi_0$.
  A different choice of constants corresponds to a gauge transformation $\psi_i \to g_i^{-1} \psi_i,$
  where $g_i$ is a diagonal matrix.
  \item We find the matrix $L'_i$ from the equation $\psi_i = (L'_i)^{-1} \psi_{i+1}$.
  This vector equation is equivalent to 3 scalar ones:
  \begin{equation}\label{find-abcde}
  \psi_{i,1} = c'_i \psi_{i+1,3}, \quad \psi_{i,2} = d'_i \psi_{i+1,1} + b'_i \psi_{i+1,3},\quad \psi_{i,3} = e'_i z \psi_{i+1,2}+a'_i \psi_{i+1,3}.
  \end{equation}
  One can check using Proposition~\ref{psi-divisor} that these equations
 determine the values $a'_i$,$b'_i$,$c'_i$,$d'_i$,$e'_i$ uniquely for each $i$.
 They do not vanish for generic spectral data.
  A gauge transformation at the previous
 step corresponds to the transformation $L'_i \to g_{i+1} L'_i g_i^{-1}$.
 \end{enumerate}
  The remaining part is to prove that $\prod_{i=0}^{n-1} c'_i d'_i e'_i = C^{-3}$ and that  a different choice of a divisor $D$ at the first step
  only changes the matrices $L'_i,\;i\ge0,$ up to gauge transformations.

  By construction, we have $\psi_n = T'_0 \psi_0 = Ck \psi_0$, i.e., $\det{(T'_0 - CkI)}=0$, where $T'_0 = L'_{n-1} L'_{n-2} ... L'_0$.
  At the same time, $\Gamma$ is the spectral curve of $T'_0$, i.e., $\det{(T'_0 - \tilde{C}kI)}=0$, where
  $\tilde{C}^{-3}=(z^n \det{T'_0(z)})^{-1} = \prod_{i=0}^{n-1} c'_i d'_i e'_i$.
  Now the required identity follows from $C=\tilde{C}$.

  Assume that we have a divisor $D'$ of degree $g+2$ equivalent to $D$. Two divisors are equivalent if and only if
  there is a meromorphic function $f$ on $\Gamma$ with zeroes at $D$ and with poles at $D'$.
  Therefore, a choice of the divisor $D'$ instead of $D$ at step 1 is equivalent to multiplying all functions
  $\psi_i, 0 \le i \le n,$ by the function $f$ at step 2. This multiplication does not change the matrices $L'_i$,
  which we obtain at step 3.
 \end{proof}

 \begin{prop}\label{general-lax}
  Any generic sequence of $n$ matrices:
  \[
  L'_j(z) =
  \begin{pmatrix}
   0 & 0 & c'_j\\
   d'_j & 0 & b'_j\\
   0 & e'_j z & a'_j
  \end{pmatrix}^{-1},
  0 \le j \le n-1, \qquad
  \prod_{i=0}^{n-1} c'_i d'_i e'_i = C^{-3}
  \]
  may be transformed to a unique sequence of matrices $L_j(z)$ (when $n \ne 3m$ and $C=1$) or $\tilde{L}_j(z)$ (for any $n$ and $C = (-1)^n J_q/I_q$)
  with help of gauge transformations: $L'_j \to g_{j+1} L'_j g_j^{-1},$ where
  $g_j = \text{diag}(\alpha_j,\beta_j,\gamma_j)$ $(0 \le j \le n-1, \; g_n=g_0)$ are diagonal matrices.
  Both $L_j(z)$ and $\tilde{L}_j(z)$ are defined in Theorem~\ref{lax-func}, and $C=(z^n \det{T'_{i,t}(z)})^{1/3}$.
  \end{prop}
 \begin{proof}
 The equation $L_j = g_{j+1} L'_j g_j^{-1}$ reads as:
 \[
 \begin{pmatrix}
  0 & 0 & 1\\
  1 & 0 & b_j\\
  0 & z & a_j
 \end{pmatrix} =
  g_{j+1}
  \begin{pmatrix}
   0 & 0 & c'_j\\
   d'_j & 0 & b'_j\\
   0 & e'_j z  & a'_j
  \end{pmatrix}
  g_j^{-1},
 \]
 and it implies a system of equations for $\alpha_j,\beta_j,\gamma_j,0 \le j \le n-1$:
 \begin{equation}\label{gauge-ab}
 c'_j \dfrac{\alpha_j}{\gamma_{j+1}}=d'_j \dfrac{\beta_j}{\alpha_{j+1}}=e'_j \dfrac{\gamma_j}{\beta_{j+1}}=1.
 \end{equation}
 Since these gauge transformations do not change the constant $C$, a necessary condition for the existence of solutions
 is $C=1$, or, equivalently, $\prod_{i=0}^{n-1} c'_i d'_i e'_i = 1.$
 One can check that it is also a sufficient condition, and the latter system of equations always has a one-parameter family
 of solutions if $n \ne 3m$. The parameter appears because a multiplication of all
 matrices $g_j$ by an arbitrary constant: $g_j \to \mu g_j$ leaves the above equations invariant.
 The variables $a_j,b_j$ are independent on $\mu$ due to their defining equations:
 $$
 a_j = a'_j \dfrac{\gamma_j}{\gamma_{j+1}}, \qquad b_j = b'_j \dfrac{\beta_j}{\gamma_{j+1}}.
 $$
 The reduction to the Lax matrix $\tilde{L}_j(z)$ may be done in a similar way:
 the equations $\tilde{L}_j = g_{j+1} L'_j g_j^{-1},\: 0 \le j \le n-1,$ are equivalent to
 a system of equations
 $$
 c'_j \dfrac{\alpha_j}{\gamma_{j+1}}=b'_j \dfrac{\beta_j}{\gamma_{j+1}}=-a'_j \dfrac{\gamma_j}{\gamma_{j+1}}, \quad e'_j \dfrac{\gamma_j}{\beta_{j+1}}=1, \quad 0 \le j \le n-1,
 $$
 which has a one-parameter family of solutions for any $n$ provided that $C = (-1)^n J_q/I_q$, or,
 equivalently, $\prod_{j=0}^{n-1} a'_j=(-1)^n \prod_{j=0}^{n-1} b'_j e'_j$. The variables $x_j,y_j$ are given by
 the formulas: $x_j = -d'_{j-2} \beta_{j-2}/\alpha_{j-1}$, $y_j = \gamma_{j-1}/(c'_{j-2} \alpha_{j-2})$.
 ``Generic'' hypothesis means that all variables $a'_j,b'_j,c'_j,d'_j,e'_j,\: 0 \le j \le n-1,$ should be non-zero.
 \end{proof}
 \begin{remark}
 If $n$ is a multiple of $3$, equations~(\ref{gauge-ab}) do not always have a solution. This is a manifestation
 of the fact that the coordinates $a_i,b_i$ work only when $n \ne 3m$.
 \end{remark}

 \subsection{Time evolution.}

 The remaining part of this section is to describe the time evolution of the pentagram map and to prove:
 \begin{repthma}{time-evol}
   Let $[D_{0,0}] \in J(\Gamma)$ be the point that corresponds to a generic twisted polygon at time $t=0$ after applying
  the spectral map, and $[D_{i,t}]$ be the point describing the twisted polygon at an integer time $t$.
  Then the equivalence class of the pole divisor $D_{i,t}$ of $\bar{\psi}_{i,t}$ is related to $[D_{0,0}]$ by the formulas:
  \begin{itemize}
   \item when $n$ is odd,
   $$
   [D_{i,t}] = [D_{0,0} - tO_1 + iO_2 + (t-i)W_2] \in J(\Gamma),
   $$
   \item when $n$ is even,
   $$
   [D_{i,t}] = \left[ D_{0,0} - tO_1 + \lfloor \dfrac{1+i}{2}\rfloor O_2 + \lfloor \dfrac{i}{2} \rfloor O_3 +
   \lfloor \dfrac{1+t-i}{2} \rfloor W_2 + \lfloor \dfrac{t-i}{2} \rfloor W_3 \right],
   $$
  \end{itemize}
  where $\deg{D_{i,t}} = g+2$, and $D_{0,0} \equiv D$ determines the point in $J(\Gamma)$ at $t=0$;
  provided that the corresponding spectral data remains generic up to time $t$.
  For odd $n$ the time evolution in $J(\Gamma)$ takes place along a straight line, whereas for even $n$ the evolution
  goes along a ``staircase'' (i.e., its square goes along a straight line).
\end{repthma}

 The time evolution of the pentagram map is described by the equation: $\psi_{i,t+1} = P_{i,t} \psi_{i,t}$,
 where $t$ is an integer parameter. The value $t=0$ corresponds to an initial $n$-gon.
 Proposition~\ref{ag-time} describes its time evolution at the level of divisors:
 \begin{prop}\label{ag-time}
  Generically, the divisors of the functions $\psi_{i,t,j}, 0 \le i \le n-1, 1 \le j \le 3$ have the following properties:
  \begin{itemize}
  \item when $n$ is odd,
  $$
  (\psi_{i,t,1}) \ge -D + tO_1 + (1-i)O_2+(1+i-t)W_2,
  $$
  $$
  (\psi_{i,t,2}) \ge -D + tO_1 - iO_2+W_1+(1+i-t)W_2,
  $$
  $$
  (\psi_{i,t,3}) \ge -D + tO_1 + (2-i)O_2+(i-t)W_2,
  $$
  \item when $n$ is even,
  $$
  (\psi_{i,t,1}) \ge -D + tO_1 + \lfloor\dfrac{1-i}{2}\rfloor O_2+\lfloor\dfrac{2-i}{2}\rfloor O_3+\lfloor\dfrac{1+i-t}{2}\rfloor W_2+\lfloor\dfrac{2+i-t}{2}\rfloor W_3,
  $$
  $$
  (\psi_{i,t,2}) \ge -D + tO_1 + \lfloor\dfrac{-i}{2}\rfloor O_2+\lfloor\dfrac{1-i}{2}\rfloor O_3+W_1+\lfloor\dfrac{1+i-t}{2}\rfloor W_2+\lfloor\dfrac{2+i-t}{2}\rfloor W_3,
  $$
  $$
  (\psi_{i,t,3}) \ge -D + tO_1 + \lfloor\dfrac{2-i}{2}\rfloor O_2+\lfloor\dfrac{3-i}{2}\rfloor O_3+\lfloor\dfrac{i-t}{2}\rfloor W_2+\lfloor\dfrac{1+i-t}{2}\rfloor W_3.
  $$
  \end{itemize}
 \end{prop}
 \begin{proof}
 Since the matrix $P_{i,t}$ is not available when $n=3m$, we use the coordinates $x_i,y_i$ and the matrix $\tilde{P}_{i,t}$ for the proof of this proposition.
 The vectors $\psi_{i,t}$ and $\tilde{\psi}_{i,t}$ are related by diagonal gauge matrices (see Remark~\ref{gauge-ab-xy}),
 therefore the vectors $\psi_{i,t}$ in the coordinates $a_i,b_i$ will have the same divisor structure when $n \ne 3m$.

 Proposition~\ref{psi-divisor} establishes the properties of the divisors of the functions $\tilde{\psi}_{i,t,j},\; 1 \le j \le 3,$ when $t=0$.
 To prove similar inequalities for $t>0$, we write out the components of the vector equation $\tilde{\psi}_{i,t+1} = \tilde{P}_{i,t} \tilde{\psi}_{i,t}$
 using an explicit formula for $\tilde{P}_{i,t}(z)$ from Theorem~\ref{lax-func} and count the orders of poles and zeroes of the components of the vector $\tilde{\psi}_{i,t+1}$.
 Note that it is sufficient to consider the cases $t=0$ and $t=1$.

 Consider, for example, the multiplicity of the function $\psi_{i,t+1,1}$ at the point $W_3$.
 Using the formula for $\tilde{P}_{i,t}$, we obtain: $\tilde{\psi}_{i,t+1,1} = (1-x_{i+2} y_{i+2})(\tilde{\psi}_{i,t,1}+\tilde{\psi}_{i,t,3})$.
 Therefore, when $t=0$ Proposition~\ref{psi-divisor} implies that $\tilde{\psi}_{i,1,1}$ has multiplicity $k$ for $i=2k$
 and $k+1$ for $i=2k+1$ at the point $W_3$. It equals $k+1$ for the function $\tilde{\psi}_{i,0,1}$ in both cases.
 This change is described by the divisor formula in the statement of the proposition.
 Other cases are treated in the same way.
 Two auxiliary formulas are used in the proof:
 $$
 \tilde{\psi}_{i,t,1} + \tilde{\psi}_{i,t,3} = O(z) \text{ at } O_1, \text{ and }
 x_{i+1} y_{i+1} \tilde{\psi}_{i,t,1} + \tilde{\psi}_{i,t,3} = O(1/z) \text{ at } W_1.
 $$
 These formulas follow from asymptotic expansions of the matrices $\tilde{T}_{i,t}$ in the same way as do similar formulas in the coordinates $a_i,b_i$.
 \end{proof}
 Propositions~\ref{ag-int},~\ref{general-lax}, and~\ref{ag-time} allow us to reconstruct the time evolution
 of an $n$-gon completely.

Now we are in a position to prove Theorem~\ref{time-evol} itself:
\begin{proof}
The vector functions $\psi_{i,t}$ with $i,t \ne 0$ are not normalized. The normalized vectors are equal to
$\bar{\psi}_{i,t} = \psi_{i,t}/f_{i,t},$ where $f_{i,t}=\sum_{j=1}^3 \psi_{i,t,j}$.
Proposition~\ref{ag-time} allows us to find the divisor of each function $f_{i,t}$:
\begin{itemize}
\item for odd $n$,
$$
 (f_{i,t}) = D_{i,t} - D_{0,0} + tO_1 - iO_2 + (i-t)W_2,
$$
\item for even $n$,
$$
 (f_{i,t}) = D_{i,t} - D_{0,0} + tO_1 + \lfloor\dfrac{-i}{2}\rfloor O_2 + \lfloor\dfrac{1-i}{2}\rfloor O_3 +
 \lfloor\dfrac{i-t}{2}\rfloor W_2 + \lfloor\dfrac{1+i-t}{2}\rfloor W_3.
$$
\end{itemize}
Since the divisor of any meromorphic function is equivalent to $[0]$, the result of the theorem follows.
\end{proof}

\begin{remark}\label{psi1}
Although the pentagram map preserves the spectral curve, it exchanges the marked points.
The ``staircase'' on the Jacobian appears after the identification of curves with different marking.
If we use a different normalization $\psi_{0,0,0} \equiv 1$ (i.e., if we divide the vector function $\psi_{0,0}$
by the first component instead of the sum of all components), the divisor $D \equiv D_{0,0}$ becomes:
\begin{itemize}
 \item $D = D_g + O_2 + W_2$ for odd $n$,
 \item $D = D_g + O_3 + W_3$ for even $n$,
\end{itemize}
where $D_g$ is a generic divisor of degree $g$ on $\Gamma$.
Note that it does not change Propositions~\ref{psi-divisor} and~\ref{ag-time}, because
only one vector $ \psi_{i,t}$ with $i=t=0$ is normalized.
\end{remark}

\section{Periodic case - closed polygons}

 In this section we prove:
 \begin{repthma}{closed-poly}
  For generic closed polygons the pentagram map is defined only for $n \ge 5$.
  Closed polygons are singled out by the condition that $(z,k)=(1,1)$ is
  a triple point of $\Gamma$. The latter is equivalent to 5 linear relations on $I_j,J_j$:
  \begin{equation}\label{lin-rel}
  \sum_{j=0}^q I_j = \sum_{j=0}^q J_j = 3, \quad \sum_{j=0}^q j I_j = \sum_{j=0}^q j J_j = 3q-n, \quad
  \sum_{j=0}^q j^2 I_j = \sum_{j=0}^q j^2 J_j.
  \end{equation}
  The genus of $\Gamma$ drops to $g=n-5$ when $n$ is even, and to $g=n-4$ when $n$ is odd.
  The dimension of the Jacobian $J(\Gamma)$ drops by $3$ for closed polygons.
  Theorem~\ref{spectral-th} holds with this genus adjustment on the space $\mathcal{C}_n$, and Theorem~\ref{time-evol} holds
  verbatim for closed polygons.
 \end{repthma}
 \begin{proof}
 The monodromy matrix from the definition of the twisted $n$-gon equals $T_{0,t}(1)$.
 Clearly, an $n$-gon is closed if and only if $T_{0,t}(1)=C \cdot Id$
 ($C=1$ for the Lax matrix $L_{i,t}$).
 The latter condition implies that $(z,k)=(1,1)$ is a self-intersection point for $\Gamma$.
 The algebraic conditions implying that $(1,1)$ is a triple point are:
 \begin{itemize}
 \item $R(1,1) = 0$,
 \item $\partial_k R(1,1) = \partial_z R(1,1) = 0$,
 \item $\partial_k^2 R(1,1) = \partial_z^2 R(1,1) = \partial_{kz}^2 R(1,1) = 0$.
 \end{itemize}
 They are equivalent to 5 linear relations among $I_j,J_j$:
 $$ \sum_{j=0}^q I_j = \sum_{j=0}^q J_j = 3, \quad \sum_{j=0}^q j I_j = \sum_{j=0}^q j J_j = 3q-n, \quad
 \sum_{j=0}^q j^2 I_j = \sum_{j=0}^q j^2 J_j.$$
 Equivalent relations were found in Theorem 4 in~\cite{ost10} (but only for the variables $a_i,b_i$).

 The proofs of Theorems~\ref{spectral-th} and~\ref{time-evol} apply, mutatis mutandis, to
 the periodic case with one change: a count of the number of branch points $\nu$ of $\Gamma$ and
 the corresponding calculation for the genus $g$ of $\Gamma$.  Generic spectral data for closed polygons
 consists of spectral functions that are singular at the point $(z,k)=(1,1)$ in addition to the points $O_i,W_i$,
 whereas the restrictions on the divisors $D \equiv D_{0,0}$ are the same as for twisted polygons.

 As before, the function $\partial_k R$ has poles of total order $3n$ above $z=0$, and zeroes of total order $n$
 about $z=\infty$. Now since $R(z,k)$ has a triple point $(1,1)$, $\partial_k R$ has a double zero at $(1,1)$.
 But $z=1$ is not a branch point of the normalization $\Gamma$. Consequently, $\partial_k R$ has double
 zeroes on 3 sheets of $\Gamma$ above $z=1$.
 The Riemann-Hurwitz formula for even $n$ becomes: $2-2g=6-\nu, \; \nu = 3n-n-6=2n-6$, and for odd $n$:
 $2-2g=6-\nu, \; \nu = 3n-n-6+2=2n-4$. Therefore, we have $g=n-5$ for even $n$, and $g=n-4$ for odd $n$.

 Remark~\ref{dim-count} implies that there are no other relations among $I_i,J_i,0 \le i \le q,$ except for~(\ref{lin-rel}) in
 the periodic case. The dimension of the Jacobian $J(\Gamma)$ is $3$ less than for twisted polygons.
 It is consistent with the fact that closed polygons form a subspace of codimension $8$ in $\mathcal{P}_n$.
 \end{proof}

 \begin{corollary}\label{cor-periodic}
 The dimension of the phase space $\mathcal{C}_n$ in the periodic case is $2n-8$.
 In the complexified case, a Zariski open subset of $\mathcal{C}_n$ is fibred over the base of dimension $2q-3$.
 The coordinates on the base are $I_j,J_j,0 \le j \le n-1,$ subject to the constraints from Theorem~\ref{closed-poly}.
 The fibres are Zariski open subsets of Jacobians (complex tori) of dimension $2q-3$ for odd $n$, and of dimension $2q-5$ for even $n$.
 Note that the restriction of the symplectic form (corresponding to the Poisson structure on the symplectic leaves)
 to the space $\mathcal{C}_n$ is always degenerate, therefore the Arnold-Liouville theorem is not directly
 applicable for closed polygons. Nevertheless, the algebraic-geometric methods guarantee
 that the pentagram map exhibits quasi-periodic motion on a Jacobian.
 \end{corollary}
 \begin{remark}
 The dimension of the tori is one when $n=5$ (for pentagons). The motion on them turns out to be periodic with period $5$.
 On the other hand, the pentagram map is known to be the identity map, see~\cite{s92}. The discrepancy appears because pentagons with a different
 numbering of vertices are considered to be the same in~\cite{s92}, but different in our paper (i.e., if we enumerate the vertices
 from $1$ to $5$ and then perform a cyclic permutation, these pentagons will not be equivalent).
 \end{remark}

 \section{The symplectic form and action-angle variables}\label{kpform-sec}

 \begin{definition}[\cite{KP97,KP98}]\label{KP-universal}
 \emph{Krichever-Phong's universal formula} defines a pre-symplectic form on the space of Lax operators, i.e.,
 on the space $\mathcal{P}_n$.
 It is given by the expression:
 $$
 \omega = -\dfrac{1}{2} \sum_{z=0,\infty} {\text{res}} \thinspace \text{Tr}\left( \Psi_0^{-1} T_0^{-1} \delta T_0 \wedge \delta \Psi_0
 \right) \dfrac{dz}{z}.
 $$
 The matrix $\Psi_{0,t}$ is defined in Proposition~\ref{deg-D}. In this section we drop the index $t$, because
 all variables correspond to the same moment of time.

 The \emph{leaves} of the 2-form $\omega$ are defined as submanifolds of $\mathcal{P}_n$, where the expression $\delta \ln{k} dz/z$
 is holomorphic. The latter expression is considered as a one-form on the spectral curve $\Gamma$.
 \end{definition}
 \begin{remark}\label{omega-gauge-inv}
 A heuristic principle justified by many examples is that when $\omega$ is restricted to these leaves,
 it becomes a symplectic form of rank $2g$, where $g$ is the genus of $\Gamma$.
 Moreover, one can prove~(\cite{IK04}) that $\omega$ does not depend on the normalization of the eigenvectors
 used to construct the matrix $\Psi_{0,t}$, and on gauge transformations
 $L_j \to g_{j+1} L_j g_j^{-1}, \; g_j \in GL(3,\C)$, when restricted to the leaves.
 \end{remark}
 \begin{remark}
 There exist different variations of the universal formula, which provide $2$ or even more compatible Hamiltonian
 structures for some integrable systems. However, it seems likely that other modifications
 of the universal formula do not lead to \emph{symplectic} forms for the pentagram map.
 \end{remark}

 \begin{repthma}{kp-form}
  Krichever-Phong's pre-symplectic 2-form  on the space $\mathcal{P}_n$ turns out to be
 a symplectic form of rank $2g$ after the restriction to the leaves:
 $\delta I_q = \delta J_q = 0$ for odd $n$, and
 $\delta I_0 = \delta I_q = \delta J_0 = \delta J_q = 0$ for even $n$.
 These leaves coincide with the symplectic leaves of the Poisson structure found in~\cite{ost10}.
 When restricted to the leaves, the 2-form $\omega$ defined above equals:
  $$
  \omega = \sum_{j=0}^{q-1} \delta \ln{x_{2j+1}} \wedge \delta \ln{\left( \prod_{k=0}^j x_{2k} \right)} -
  \sum_{j=0}^{q-1} \delta \ln{y_{2j+1}} \wedge \delta \ln{\left( \prod_{k=0}^j y_{2k} \right)}.
  $$
 This symplectic form is invariant under the pentagram map and coincides
 with the inverse of the Poisson structure restricted to the symplectic leaves.
 It has natural Darboux coordinates, which turn out to be action-angle coordinates for the pentagram map.
 \end{repthma}

\begin{proof}
 In this proof we again invoke Remark~\ref{gauge-ab-xy}:
 the formula for $\omega$ is algebraic, and our proof never uses the ``non-divisibility by $3$'' condition.
 Therefore, our computation of $\omega$ in the coordinates $a_i,b_i$ is formally valid for all $n$.
 This remark also implies that $\tilde{T}_0 = (J_q/I_q) g_0^{-1} T_0 g_0$.
 Remark~\ref{omega-gauge-inv} and the fact that $\delta (J_q/I_q)=0$ imply that our formal computation gives
 a correct symplectic structure for all $n$ when $\omega$ is written in the coordinates $x_i,y_i$.

 First we find the equations that define the leaves of the 2-form $\omega$.
 \begin{lemma}\label{s-leaves}
 The one-form $\delta \ln{k} dz/z$ is holomorphic on the spectral curve $\Gamma$ when
 restricted to the leaves: $\delta I_q = \delta J_q = 0$ for odd $n$, and
  $\delta I_0 = \delta I_q = \delta J_0 = \delta J_q = 0$ for even $n$.
 \end{lemma}
 \begin{proof}
  These conditions follow immediately from the definition of the leaves and Lemma~\ref{spec-sing}.
 For example, at the point $O_1$ we have
 $$
 \delta \ln{k_1} \dfrac{dz}{z} = \left(\dfrac{1}{z} \dfrac{\delta I_q}{I_q} + O(1) \right) dz.
 $$
 Clearly, this one-form is holomorphic in $z$ if and only if $\delta I_q = 0$.
 Similarly, we obtain $\delta I_q=0$ at the point $O_2$ for odd $n$. One has to keep in mind that the local
 parameter is $\lambda^{1/2}$ there.
 \end{proof}
Now we introduce the action-angle coordinates. Their construction is universal, see, in particular,
the proof of Corollary 4.2 in~\cite{IK02}.
 \begin{lemma}
  The rank of $\omega$ is $2g$ when it is restricted to the leaves of Lemma~\ref{s-leaves}.
 $$
 \omega=\sum_{i=1}^{g} \delta \mathbf{I}_i \wedge \delta\varphi_i, \text{ where }
 \mathbf{I}_i = \oint_{a_i} \ln{k} dz/z, \quad \varphi_i = \sum_{s=1}^{g+2} \int^{\gamma_s} d\omega_i,
 $$
 the points $\gamma_s \in \Gamma,\; 1 \le s \le g+2,$ are the points of the divisor $D_{0,0}$ from Proposition~\ref{deg-D},
 and $\mathbf{I}_i,\varphi_i$ are called action-angle coordinates (the angle variables $\varphi_i,\; i\ge 0,$ are defined on the Jacobian $J(\Gamma)$).
 \end{lemma}
 \begin{proof}
 Since the one-form $\delta \ln{k} dz/z$ is holomorphic on $\Gamma$, it can be represented as a sum of the basis
 holomorphic differentials:
\begin{equation}\label{kpformtmp4}
 \delta \ln{k} dz/z = \sum_{i=1}^{g} \delta \mathbf{I}_i d\omega_i,
 \end{equation}
 where $g$ is the genus of $\Gamma$.
The coefficients $\mathbf{I}_i$ may be found by integrating the last expression over the basis cycles $a_i$ of $H_1(\Gamma)$:
$$
\mathbf{I}_i = \oint_{a_i} \ln{k} dz/z.
$$
According to formula~(5.7) in~\cite{KP00}, we have:
$$
\omega=\sum_{i=1}^{g+3} \delta \ln{k(\gamma_i)} \wedge \delta \ln{z(\gamma_i)},
$$
where the points $\gamma_i \in \Gamma,\; 1 \le i \le g+3,$ constitute the pole divisor $D_{0,0}$ of the
normalized Floquet-Bloch solution $\psi_{0,0}$ from Proposition~\ref{deg-D}.

After a rearrangement of terms, we obtain:
$$
\omega = \delta\left( \sum_{s=1}^{g+2} \int^{\gamma_s} \delta \ln{k} dz/z \right)=
\delta\left( \sum_{s,i} \int^{\gamma_s} \delta \mathbf{I}_i d\omega_i \right) = \sum_{i=1}^{g} \delta \mathbf{I}_i \wedge \delta\varphi_i,
$$
where
$$
\varphi_i = \sum_{s=1}^{g+3} \int^{\gamma_s} d\omega_i
$$
are coordinates on the Jacobian of $J(\hat{\Gamma})$.
The variables $\mathbf{I}_i$ and $\varphi_i$ are known as action-angle variables.

Let us show that the functions $\mathbf{I}_i$ are independent.
If they are not, then there exists a vector $v$ on the space ${\mathcal P}_n$, such that $\delta \mathbf{I}_i(v)=0$ for all $i$.
Then it follows from~(\ref{kpformtmp4}), that $\partial_v k \equiv 0$.
If we apply $\partial_v$ to $R(z,k)$, we conclude that $k$ satisfies an algebraic equation
of degree $2$, which is impossible, since $\Gamma$ is a $3$-sheeted cover of $z$-plane.
\end{proof}

Finally, we proceed to the computation of $\omega$. Note that
$$
\text{Tr}\left( \Psi_0^{-1} T_0^{-1} \delta T_0 \wedge \delta \Psi_0 \right) =
\sum_{k=0}^{n-1} \text{Tr}\left( \Psi_0^{-1} L_0^{-1}...L_k^{-1} \delta L_k L_{k-1}...L_0 \wedge \delta \Psi_0 \right)=
$$
$$
= \sum_{k=0}^{n-1} \text{Tr}\left( \Psi_k^{-1} L_k^{-1} \delta L_k \wedge \delta \Psi_k \right) -
\sum_{k=0}^{n-2} \text{Tr}\left( L_0^{-1}...L_k^{-1} \delta L_k \wedge \delta(L_{k-1}...L_0) \right),
$$
where $\Psi_k = L_{k-1}...L_0 \Psi_0$ (this transformation is similar to the one used in~\cite{IK04}).
Notice that the last sum does not have any poles except at the points $z=0$ and $z=\infty$ and vanishes after
the summation over both residues. Therefore,
$$
\omega = -\dfrac{1}{2} \sum_{j=0}^{n-1} \underset{0,\infty}{\text{res}} \thinspace \text{Tr}\left( \Psi_j^{-1} L_j^{-1} \delta L_j \wedge \delta \Psi_j \right) \dfrac{dz}{z}.
$$
To compute $\omega$, we use a normalization of $\psi_0$ in which $\psi_{0,1} \equiv 1$. It corresponds to the case when
the first line of $\Psi_0$ is $(1,1,1)$. The matrices $\Psi_j, j > 0,$ are not normalized.
A normalized matrix $\bar{\Psi}_j, j > 0,$ is related to
$\Psi_j$ by a diagonal matrix $F_j$: $\bar{\Psi}_j = \Psi_j F_j$.
The matrices $F_j,j>0,$ may have poles or zeroes at $z=0,\infty$. We have the formula:
$$
\text{Tr}\left( \Psi_j^{-1} L_j^{-1} \delta L_j \wedge \delta \Psi_j \right) =
\text{Tr}\left( \bar{\Psi}_j^{-1} L_j^{-1} \delta L_j \wedge \delta \bar{\Psi}_j \right) -
\text{Tr}\left( \bar{\Psi}_j^{-1} L_j^{-1} \delta L_j \bar{\Psi}_j \wedge \delta \ln{F_j} \right)
$$
Notice that the product $L_j^{-1} \delta L_j$ is
$$
L_j^{-1} \delta L_j =
\begin{pmatrix}
            0 & 0 & 0\\
  -\delta b_j & 0 & 0\\
  -\delta a_j & 0 & 0
\end{pmatrix},
$$
and the first line of $\delta \bar{\Psi}_j$ is always zero due to the normalization.
Consequently, we obtain the formula:
$$
\omega = \dfrac{1}{2} \sum_{j=0}^{n-1} \underset{0,\infty}{\text{res}} \thinspace
\text{Tr}\left( \bar{\Psi}_j^{-1} L_j^{-1} \delta L_j \bar{\Psi}_j \wedge \delta \ln{F_j} \right) \dfrac{dz}{z}.
$$
We can rewrite the last formula as:
$$
\omega = \dfrac{1}{2} \sum_{j=0}^{n-1} \sum_i \underset{O_i,W_i}{\text{res}} \thinspace
\text{Tr}\left( \psi^*_j L_j^{-1} \delta L_j \bar{\psi}_j \wedge \delta \ln{f_j} \right) \dfrac{dz}{z},
$$
where $\psi^*_j$ is an eigen-covector: $\psi^*_j T_j = k \psi^*_j$. Covectors are normalized by $\psi^*_j \bar{\psi}_j = 1$,
and $\bar{\psi}_{j,1} \equiv 1$.
One can check that $\psi^*_j L_j^{-1} \delta L_j \bar{\psi}_j = -\psi^*_{j,3} \delta a_j - \psi^*_{j,2} \delta b_j$.
The formula for $\omega$ becomes:
\begin{equation}\label{main-sf}
\omega = -\dfrac{1}{2} \sum_{j=0}^{n-1} \sum_i \underset{O_i,W_i}{\text{res}} \thinspace
(\psi^*_{j,3} \delta a_j + \psi^*_{j,2} \delta b_j) \wedge \delta \ln{f_j} \dfrac{dz}{z} =
\sum_i \omega_{O_i} + \sum_i \omega_{W_i}.
\end{equation}
We use formula~(\ref{main-sf}) to compute $\omega$.
We compute the terms $\omega_{O_i}$ and $\omega_{W_i}$ with different $i$ separately in Lemmas~\ref{omega-O1}-\ref{omega-W23}.
One can show that their sum equals:
$$
 \omega = \sum_{(i,j) \in \Lambda} (\delta \ln{a_i} \wedge \delta \ln{a_j} - \delta \ln{b_i} \wedge \delta \ln{b_j}),
$$
where the set $\Lambda$ consists of pairs $(i,j), 0 \le i \le n-1, i < j \le n-1$, such that
either both $i$ and $j$ are even, or $i$ is odd and $j$ is arbitrary.
The integrals of motion for the Lax matrix $\tilde{L}_{i,t}(z)$ are related
to the Casimirs $E_n,O_n,E_{n/2},O_{n/2}$ found in Corollary 2.13 in~\cite{ost10} in the following way:
$$
\text{for any $n$,}\quad E_n = (-1)^n \dfrac{J_q}{I_q^2}, \quad O_n = \dfrac{I_q}{J_q^2};\quad \text{for even $n$,} \quad
 E_{n/2} = (-1)^q \dfrac{I_0}{I_q}, \quad O_{n/2} = (-1)^q \dfrac{J_0}{J_q}.
$$
Clearly, these Casimirs define the same symplectic leaves as Lemma~\ref{s-leaves}.
One can show using formulas~(\ref{def-xy}) that on these leaves $\omega$ equals
$$
\omega_0 = \sum_{j=0}^{q-1} \delta \ln{x_{2j+1}} \wedge \delta \ln{\left( \prod_{k=0}^j x_{2k} \right)} -
\sum_{j=0}^{q-1} \delta \ln{y_{2j+1}} \wedge \delta \ln{\left( \prod_{k=0}^j y_{2k} \right)}.
$$
On the leaves, its inverse equals the Poisson brackets (2.16) in~\cite{ost10}:
$$
\{ x_i, x_j \} = (\delta_{i,j-1} - \delta_{i,j+1}) x_i x_j, \qquad
\{ y_i, y_j \} = (\delta_{i,j+1} - \delta_{i,j-1}) y_i y_j,
$$
and all other brackets vanish.
Note that since the symplectic leaves for these Poisson brackets have positive codimension,
the corresponding $2$-form is not unique, and $\omega_0$ represents one of the possible $2$-forms.
\end{proof}

\section{Appendix}
In this appendix we prove Lemmas~\ref{T0-asymp}-\ref{omega-W23}, which complete the proof of Proposition~\ref{psi-divisor} and Theorem~\ref{kp-form}.
\begin{lemma}\label{T0-asymp}
 When $n=2q$, the expansion of $T_{0,t}(z)$ at $z=0$ is:
 \[
 T_{0,t}(z) =
 \begin{pmatrix}
  (-1)^q \prod_{i=0}^{q-1} a_{2i}   & 0 & (-1)^{q-1} \prod_{i=1}^{q-1} a_{2i}\\
  O(1) & (-1)^q \prod_{i=0}^{q-1} a_{2i+1} & O(1)\\
    0    & 0 & 0
 \end{pmatrix}
 \dfrac{1}{z^q} + O\left( \dfrac{1}{z^{q-1}} \right),
 \]
 and the expansion of $T^{-1}_{0,t}(z)$ at $z=\infty$ is:
 \[
 T^{-1}_{0,t}(z) =
 \begin{pmatrix}
  0                          & \prod_{i=1}^{q-1} b_{2i} & 0\\
  0                          & \prod_{i=0}^{q-1} b_{2i} & 0\\
  \prod_{i=1}^{q-1} b_{2i-1} & O(1) & \prod_{i=1}^q b_{2i-1}
 \end{pmatrix}
 z^q + O(z^{q-1}).
 \]
 When $n=2q+1$, we have:
  \[
 T_{0,t}(z) =
 \begin{pmatrix}
    O(z^{-q}) & \dfrac{\prod_{i=1}^q (-a_{2i-1})}{z^q} + O(z^{1-q}) & O(z^{-q})\\
    \dfrac{\prod_{i=0}^q (-a_{2i})}{z^{q+1}} + O(z^{-q}) & O(z^{-q}) & \dfrac{\prod_{i=1}^q (-a_{2i})}{z^{q+1}} + O(z^{-q})\\
    O(z^{-q})    & O(z^{1-q}) & O(z^{-q})
 \end{pmatrix},
 \]
 \[
 T^{-1}_{0,t}(z) =
 \begin{pmatrix}
  O(z^q) & O(z^q) & z^q \prod_{i=1}^q b_{2i} + O(z^{q-1})\\
  z^q \prod_{i=0}^{q-1} b_{2i} + O(z^{q-1}) & O(z^q) & z^q \prod_{i=0}^q b_{2i} + O(z^{q-1})\\
  O(z^q) & z^{q+1} \prod_{i=1}^q b_{2i-1} + O(z^q) & O(z^q)
 \end{pmatrix}.
 \]
\end{lemma}
\begin{proof}
Let us prove the first formula for $n=2q$ (the others are proved similarly).
One can check that
 \[
 L_{j+1}(z) L_j(z) =
 \begin{pmatrix}
 -a_j & 0 & 1\\
 1+a_{j+1} b_j & -a_{j+1} & 0\\
 0 & 0 & 0
 \end{pmatrix}
 \dfrac{1}{z} + O(1) \text{ at } z=0,
 \]
 and the expansion for $T_{0,t}(z)$ follows from it.
\end{proof}

\begin{lemma}\label{omega-O1}
 The contribution from the point $O_1$ is independent on the parity of $n$ and is given by:
 $$
 \omega_{O_1} = -\dfrac{1}{2} \sum_{j=2}^{n-1} \delta \ln{a_j} \wedge \delta \ln{\left( \prod_{k=1}^j a_k \right)}.
 $$
\end{lemma}
\begin{proof}
First, we prove 2 formulas:
\begin{equation}\label{psi-O1}
\bar{\psi}_j(O_1) = \left( 1, \dfrac{1}{a_{j+1}} + b_j, a_j \right)^T, \qquad \psi^*_j(O_1)=(0,0,1/a_j),
\end{equation}
then we find $f_j(O_1)$, and compute $\omega_{O_1}$ using formula~(\ref{main-sf}).

The vectors $\psi_0$, $\psi_0^*$ and the matrix $T_0$ are related to $\bar{\psi}_j$, $\psi_j^*$, $T_j$ by
a permutation of the variables $a_0,...,a_{n-1}$ and $b_0,...,b_{n-1}$. Therefore, formulas~(\ref{psi-O1})
are equivalent to 2 formulas (which we prove below):
$$
\psi_0(O_1) = \left( 1, \dfrac{1}{a_1} + b_0, a_0 \right)^T, \quad \psi^*_0(O_1)=(0,0,1/a_0).
$$
Proposition~\ref{psi-divisor} and formulas~(\ref{findab0}) imply that $\psi_0(O_1) = (1,x,a_0)^T$ for some constant $x$.
Using the value of $T_0^{-1}$ at $z=0$:
$$
T_0^{-1}(0) =
\begin{pmatrix}
  0 & 0 & I_q/a_0\\
  0 & 0 & (1+a_1 b_0) I_q/(a_0 a_1)\\
  0 & 0 & I_q
\end{pmatrix},
$$
and the formula $T_0^{-1}(0) \psi_0(O_1) = I_q \psi_0(O_1)$, we find that $x=(1/a_1)+b_0$.

One can check that the equation $\psi^*_0 T_0 = k \psi^*_0$ implies that $\psi^*_0(O_1)=(0,0,y)$
for some constant $y$. Since $\psi^*_0 \psi_0 = 1$, we find that $y=1/a_0$.

To find $f_j(O_1)$, we have to compare $\bar{\psi}_j$ and $L_{j-1}...L_0 \psi_0$ at the point $O_1$.
One can check that $L_0 \psi_0(O_1) = (1/a_1, *, 1)^T$. Therefore, $f_1(O_1) = a_1$.
When $i>0$, we have $L_i \bar{\psi}_i = (1/a_{i+1}, *, 1)^T$. Consequently, we find that $f_i(O_1)/f_{i-1}(O_1)=a_i$.
Multiplying the latter equations with $2 \le i \le j$ by each other, we obtain that $f_j(O_1) = \prod_{k=1}^j a_k$.

Substituting $f_j(O_1)$ and $\psi^*_j(O_1)$ into formula~(\ref{main-sf}), we obtain $\omega_{O_1}$.
\end{proof}
Similarly, the contribution from the point $W_1$ is given by:
\begin{lemma}
For both even and odd $n$,
$$
\omega_{W_1} = \dfrac{1}{2} \sum_{j=1}^{n-1} \delta \ln{b_j} \wedge \delta \ln{\left( \prod_{k=0}^{j-1} b_k \right)}.
$$
\end{lemma}
\begin{proof}
In the same way as in Lemma~\ref{omega-O1}, we find that
$$
\bar{\psi}_j(W_1) = (1,0,-1/b_{j-1})^T, \qquad \psi^*_j(W_1) = (1,-1/b_j,0), \qquad f_j(W_1) = (-1)^j \prod_{k=0}^{j-1} b_k^{-1},
$$
which implies the formula for $\omega_{W_1}$.
\end{proof}
The computation at the points $O_2,O_3,W_2,W_3$ is trickier, because it differs for even and odd $n$.
\begin{lemma}\label{omega-O2}
If $n$ is odd, then
$$
\omega_{O_2} =
-\dfrac{1}{2} \sum_{j=1}^{n-1} \delta \ln{a_j} \wedge \delta \ln{\left( \prod_{k=0}^{j-1} \prod_{i=0}^q a_{k+2i} \right)}.
$$
\end{lemma}
\begin{proof}
First, we need to prove 2 formulas:
\begin{equation}\label{O2-eq1}
\bar{\psi}_j = \left( 1, \alpha_j \dfrac{1}{\sqrt{z}} + O(1),
\beta_j \sqrt{z} + O(z) \right)^T \text{ at } O_2,
\end{equation}
$$
\text{ where }\alpha_j = \dfrac{(-1)^{q+1} \prod_{i=0}^q a_{j+2i}}{\sqrt{-I_q}}, \qquad
\beta_j = \dfrac{(-1)^{q} \prod_{i=0}^{q-1} a_{j+2i}}{\sqrt{-I_q}}.
$$
\begin{equation}\label{O2-eq2}
\psi^*_j(O_2) = \left( \dfrac{1}{2},0,-\dfrac{1}{2a_j} \right).
\end{equation}
Note that a cyclic permutation of the variables $a_j \to a_{j+1}, b_j \to b_{j+1}$ (for all $j$)
permutes the eigenvectors and covectors as follows: $\bar{\psi}_j \to \bar{\psi}_{j+1}$,
$\psi^*_j \to \psi^*_{j+1}$. Therefore, we only need to find $\bar{\psi}_0$ at $O_2$ and $\psi^*_0(O_2)$
to prove formulas~(\ref{O2-eq1}) and~(\ref{O2-eq2}).

Proposition~\ref{psi-divisor} implies that $\psi_0 = (1,\alpha_0/\sqrt{z}+O(1),\beta_0 \sqrt{z}+O(z))^T$
around the point $O_2$.
Since $T_0 \psi_0 = \left(\sqrt{-I_q} z^{-n/2} + O(z^{-q})\right)\psi_0$, we find that
$$
\alpha_0 = \dfrac{(-1)^{q+1} \prod_{i=0}^q a_{2i}}{\sqrt{-I_q}}.
$$
One can check that $\left(T_0^{-1}\right)_{32} = I_q z/a_{n-1} + O(z^2)$, and
since $T_0^{-1} \psi_0 = \psi_0 O(z^{n/2})$ in the neighborhood of $O_2$, we deduce that $\beta_0 = -\alpha_0/a_{n-1}$.
Formula~(\ref{O2-eq1}) with $j=0$ is proven.

The equation $\psi^*_0 \psi_0 = 1$ implies that
$\psi^*_0 = (\alpha'+O(\sqrt{z}), \beta' \sqrt{z} + O(z), \gamma'+O(\sqrt{z}))$ at the point $O_2$.
Using the identity $\psi^*_0 T_0 = \left(\sqrt{-I_q} z^{-n/2} + O(z^{-q})\right)\psi^*_0$,
we find that
$$
 \beta' \prod_{i=0}^q (-a_{2i}) = \alpha' \sqrt{-I_q}, \qquad
 \beta' \prod_{i=1}^q (-a_{2i}) = \gamma' \sqrt{-I_q}, \qquad
 \alpha' + \beta' \dfrac{(-1)^{q+1} \prod_{i=0}^q a_{2i}}{\sqrt{-I_q}} = 1.
$$
Solving these equations for $\alpha',\beta',\gamma'$, we obtain that $\psi^*_0(O_2) = (1/2,0,-1/(2a_0))$.

Now we find the value of $\delta \ln{f_j(O_2)}$.
Since $(L_0 \psi_0)_1 = \psi_{0,2}-b_0 \psi_{0,1}$, we obtain that $\delta \ln{f_1}(O_2) = - \delta \ln{\alpha_0}$.
The argument similar to the one used in the proof of Lemma~\ref{omega-O1} implies that
$$
f_j(z) = \dfrac{z^{j/2}}{\prod_{k=0}^{j-1} \alpha_k} + O\left( z^{(j+1)/2)} \right) \text{ at } O_2.
$$
Using the condition $\delta I_q=0$, we obtain that
$$
 \delta \ln{f_j(O_2)} = -\delta \ln{\left( \prod_{k=0}^{j-1} \prod_{i=0}^q a_{k+2i} \right)}.
$$
Finally, using formula~(\ref{main-sf}), we deduce
$$
\omega_{O_2} = \dfrac{1}{2} \sum_{j=1}^{n-1} 2 \cdot \dfrac{1}{2} \delta \ln{a_j} \wedge \delta \ln{f_j(O_2)} =
-\dfrac{1}{2} \sum_{j=1}^{n-1} \delta \ln{a_j} \wedge \delta \ln{\left( \prod_{k=0}^{j-1} \prod_{i=0}^q a_{k+2i} \right)}.
$$
The coefficient ``$2$'' in the last formula appears because $O_2$ is a branch point. The local parameter around the point $O_2$
is $\sqrt{z}$, and one has to use the formula $2(d\sqrt{z})/\sqrt{z}$ instead of $dz/z$ to compute the residue at $O_2$.
\end{proof}

\begin{lemma}
If $n$ is odd, then
$$
\omega_{W_2} =
-\dfrac{1}{2} \sum_{j=1}^{n-1} \delta \ln{b_j} \wedge \delta \ln{\left(\prod_{k=0}^{j-1} \prod_{i=1}^q b_{k+2i+1} \right)}.
$$
\end{lemma}
\begin{proof}
The computation of $\omega_{W_2}$ is very similar to that of $\omega_{O_2}$ in Lemma~\ref{omega-O2}.
We compute $\psi_0$, $\psi_0^*(W_2)$, $f_1(z)$, and find the expressions for $f_j(z)$
and $\psi^*_j(W_2)$ with arbitrary $j$:
\begin{equation}\label{psi-W2}
f_j(z) = \dfrac{z^{j/2}}{\prod_{k=0}^{j-1} \beta_k} + O\left( z^{(j-1)/2)} \right) \text{ at } W_2, \qquad
\psi^*_j(W_2)=\left(0,\dfrac{1}{2b_j},0\right).
\end{equation}

From Proposition~\ref{psi-divisor} and formula~(\ref{findab0}) it follows that
$\psi_0 = (1,b_0+\beta_0/\sqrt{z}+O(1/z),\alpha_0 \sqrt{z}+O(1))^T$ near the point $W_2$,
where
$$
\beta_j = \left( \prod_{i=1}^q b_{j+2i+1} \right)/\sqrt{-J_q}
$$
From the identity $(T_0^{-1} \psi_0)_1 = k^{-1} \psi_{0,1}$ we find that
$\alpha_0 \prod_{i=1}^q b_{2i} = \sqrt{-J_q}.$
The identity $(T_0 \psi_0)_1 = k \psi_{0,1}$, along with the formulas:
$$
T_0(W_2) =
\begin{pmatrix}
  J_q & -J_q/b_0 & 0\\
  0 & 0 & 0\\
  -J_q/b_{n-1} & J_q/(b_0 b_{n-1}) & 0
\end{pmatrix}, \qquad
T_0(z)_{13} = J_q/(b_0 b_1 z) + O(z^{-2}) \text{ near } W_2,
$$
implies that $\beta_0 b_1 = \alpha_0$. Solving the above equations for $\beta_0$, we find that
$\beta_0 = \left( \prod_{i=1}^q b_{2i+1} \right)/\sqrt{-J_q}$, and the corresponding formulas for $\alpha_j, \beta_j$, and $\bar{\psi}_j$ follow.

Since $(L_0 \psi_0)_1 = \beta/\sqrt{z} + O(1/z)$, we obtain that
$$
f_1(z) = \dfrac{\sqrt{z}}{\beta_0} + O(1),\; \text{ and }\;
f_j(z) = \dfrac{z^{j/2}}{\prod_{k=0}^{j-1} \beta_k} + O\left( z^{(j-1)/2)} \right) \text{ at } W_2.
$$
Consequently, we have $\delta \ln{f_j}(W_2) = - \delta \ln{\left(\prod_{k=0}^{j-1} \beta_k \right)}$, and since
$\delta J_q = 0$ on the symplectic leaves, the formula for $\delta \ln{f_j}(W_2)$ follows.

Now we find the covector $\psi^*_0$ at the point $W_2$. The identity $\psi^*_0 \psi_0 \equiv 1$
implies that
$\psi^*_0 = (A + O(1/\sqrt{z}), B + O(1/\sqrt{z}), C/\sqrt{z} + O(1/z))$,
and that $A + b_0 B + \alpha C = 1$.
The identity $(\psi^*_0 T_0^{-1})_2 = k^{-1} \psi^*_{0,2}$ implies that $C \prod_{i=1}^q b_{2i-1} = B \sqrt{-J_q}$.
One can check that since the product $\psi^*_0 T_0$ has zero of order $n$ at $W_2$, it must be that $A=0$.
Solving the above equations for $B$, we find that $B=1/(2b_0)$, and that $\psi^*_0(W_2)=(0,1/(2b_0),0)$.

The arguments identical to those used in Lemmas~\ref{omega-O1},~\ref{omega-O2} prove formulas~(\ref{psi-W2}).
Substituting formulas~(\ref{psi-W2}) into formula~(\ref{main-sf}), we obtain:
$$
\omega_{W_2} = \dfrac{1}{2} \sum_{j=1}^{n-1} 2 \cdot \dfrac{1}{2} \delta \ln{b_j} \wedge \delta \ln{f_j}(W_2) =
-\dfrac{1}{2} \sum_{j=1}^{n-1} \delta \ln{b_j} \wedge \delta \ln{\left(\prod_{k=0}^{j-1} \prod_{i=1}^q b_{k+2i+1} \right)}.
$$
The coefficient ``-2'' appears in the last formula because
the local parameter at $W_2$ is $z^{-1/2}$, and the formula $-2 d(z^{-1/2})/z^{-1/2}$ should be used instead of $dz/z$
to compute the residue.
\end{proof}
Now we find the contribution to $\omega$ from the points $O_2,O_3,W_2,W_3$ for even $n$.
\begin{lemma}\label{omega-O23}
If $n$ is even, then
$$
\omega_{O_2} = -\dfrac{1}{2} \sum_{j=1}^{q-1} \delta \ln{a_{2j}} \wedge \delta \ln{\prod_{k=0}^{j-1} a_{2k}}, \qquad
\omega_{O_3} = -\dfrac{1}{2} \sum_{j=1}^{q-1} \delta \ln{a_{2j+1}} \wedge \delta \ln{\prod_{k=0}^{j-1} a_{2k+1}}.
$$
\end{lemma}
\begin{proof}
The substitution of the following formulas into~(\ref{main-sf}) proves the lemma:
$$
\psi^*_{2j}(O_2)=(1,0,-1/a_{2j}), \quad \psi^*_{2j+1}(O_2)=(0,0,0),
$$
$$
f_{2j}(z)= \dfrac{(-1)^j}{\prod_{k=0}^{j-1} a_{2k}} z^j + O(z^{j+1}) \text { at } O_2;
$$
$$
\psi^*_{2j+1}(O_3)=(1,0,-1/a_{2j+1}), \quad \psi^*_{2j}(O_3)=(0,0,0), \quad f_1(z) = \dfrac{z}{\eta}+O(z^2)\text{ at } O_3,
$$
$$
f_{2j+1}(z)= \dfrac{(-1)^j}{\eta \prod_{k=0}^{j-1} a_{2k+1}} z^j + O(z^{j+1}) \text{ at } O_3.
$$
Note that the parameter $\eta$ vanishes from the final formulas on the symplectic leaves.

A cyclic permutation of the variables $a_j \to a_{j+1}, b_j \to b_{j+1}$ (for all $j$) permutes the eigenvectors
and covectors as follows: $\bar{\psi}_{2j}(O_2) \to \bar{\psi}_{2j+1}(O_3)$,
$\bar{\psi}_{2j}(O_3) \to \bar{\psi}_{2j+1}(O_2)$,
$\psi^*_{2j}(O_2) \to \psi^*_{2j+1}(O_3)$, $\psi^*_{2j}(O_3) \to \psi^*_{2j+1}(O_2)$.
These permutations imply that the formulas above are equivalent to:
$$
(L_1 L_0 \psi_0)_1 = -\dfrac{a_0}{z} + O(1) \text{ at } O_2,
\quad \psi^*_0(O_2)=(1,0,-1/a_0), \quad \psi^*_0(O_3) = (0,0,0).
$$
Proposition~\ref{psi-divisor} implies that $\psi_0=(1,O(1),O(z))^T$ at the point $O_2$.
One can check that the principal part of $(L_1 L_0 \psi_0)_1$ at $O_2$ is $-a_0/z$, which implies
that $f_{j+2}(z)/f_j(z) = -z/a_j + O(z^2)$ at $O_2$ for even $j$.

Let the covector $\psi^*_0(O_2)$ be $(\alpha,\beta,\gamma)$.
The equation $(\psi^*_0 T_0)_1 = k \psi^*_{0,1}$ implies that $\beta=0$.
Since $\psi^*_0(O_2) \psi_0(O_2)=1$, we find that $\alpha=1$.
One can check that since the product $\psi^*_0 T_0^{-1}$ has zero of order $q$ at $O_2$, it must be that $\gamma = -1/a_0$.
Therefore, we obtain that $\psi^*_0(O_2)=(1,0,-1/a_0)$.

Proposition~\ref{psi-divisor} implies that $\psi_0 = (1, \eta/z+O(1),O(1))^T$ at the point $O_3$.
The principal part of $(L_0 \psi_0)_1$ at $O_3$ is $\eta/z$, and the formula for $f_1(z)$ at $O_3$ follows.
Since the product $\psi^*_0 \psi_0$ is holomorphic at $O_3$, it must be that $\psi^*_{0,2}(O_3) = 0$ and
$\psi^*_0(O_3)=(\alpha,0,\beta)$ for some $\alpha,\beta$.
One can check that the equation $\psi^*_0 T_0 = k \psi^*_0$ implies that $\alpha=\beta=0$,
thus $\psi^*_0(O_3) = (0,0,0)$.
\end{proof}

\begin{lemma}\label{omega-W23}
If $n$ is even, then
$$
\omega_{W_2} = \dfrac{1}{2} \sum_{j=1}^{q-1} \delta \ln{b_{2j}} \wedge \delta \ln{\prod_{k=0}^j b_{2k}}, \qquad
\omega_{W_3} = \dfrac{1}{2} \sum_{j=1}^{q-1} \delta \ln{b_{2j+1}} \wedge \delta \ln{\prod_{k=0}^j b_{2k+1}}.
$$
\end{lemma}
\begin{proof}
The proof of this lemma is very similar to the proof of Lemma~\ref{omega-O23}.
We prove that:
$$
\psi^*_{2j}(W_2) = (0, 1/b_{2j}, 0), \quad \psi^*_{2j+1}(W_2)=(0,0,0), \quad
f_{2j}(z) = \left( \prod_{k=1}^j b_{2k} \right) z^j + O(z^{j-1}) \text { at } W_2;
$$

$$
\psi^*_{2j+1}(W_3) = (0, 1/b_{2j+1}, 0), \quad \psi^*_{2j}(W_3)=(0,0,0),
$$
$$
f_1(W_3) = \xi, \quad
f_{2j+1}(z) = \xi \left( \prod_{k=1}^j b_{2k+1} \right) z^j + O(z^{j-1}) \text { at } W_3;
$$
and substitute these formulas into~(\ref{main-sf}).

The parameter $\xi$ vanishes from the formulas for $\omega_{W_2},\;\omega_{W_3}$ on the symplectic leaves.

A cyclic permutation $a_j \to a_{j+1}, b_j \to b_{j+1}$ (for all $j$) acts on the eigenvectors and covectors as follows:
$\bar{\psi}_{2j}(W_2) \to \bar{\psi}_{2j+1}(W_3)$,
$\bar{\psi}_{2j}(W_3) \to \bar{\psi}_{2j+1}(W_2)$,
$\psi^*_{2j}(W_2) \to \psi^*_{2j+1}(W_3)$, $\psi^*_{2j}(W_3) \to \psi^*_{2j+1}(W_2)$,
therefore we only need to prove the following:
$$
(L_0^{-1} L_1^{-1} \bar{\psi}_2)_1 = b_2 z + O(1) \text{ at } W_2,  \quad \psi^*_0(W_2)=(0,1/b_0,0), \quad \psi^*_0(W_3) = (0,0,0).
$$
Proposition~\ref{psi-divisor} implies that $\psi_0=(1,b_0+O(1/z),O(1))^T$ at the point $W_2$.
One can check that the principal part of $(L_0^{-1} L_1^{-1} \bar{\psi}_2)_1$ at $W_2$ is $b_2 z$,
which implies that $f_{j+2}(z)/f_j(z) = b_{j+2} z + O(1)$ at $W_2$ for even $j$.

Let the covector $\psi^*_0(W_2)$ be $(\alpha,\beta,\gamma)$.
One can check that the highest order terms of the equation $\psi^*_0 T_0^{-1} = k^{-1} \psi^*_0$
imply that $\alpha=\gamma=0$. Since $\psi^*_0 \psi_0 = 1$, we find that $\beta = 1/b_0$, and
$\psi^*_0(W_2)=(0,1/b_0,0)$.

Proposition~\ref{psi-divisor} implies that $\psi_0=(1,O(1),O(z))^T$ at the point $W_3$.
Therefore, $(L_0 \psi_0)_1$ is $O(1)$ at $W_3$, and we define $\xi = 1/(L_0 \psi_0)_1(W_3)$. Hence, $f_1(W_3) = \xi$.
Since $\psi^*_0 \psi_0$ is holomorphic at $W_3$, it must be that $\psi^*_0(W_3)=(\alpha,\beta,0)$ for some $\alpha,\beta$.
One can check that $\psi^*_0 T_0 = k \psi^*_0$ implies $\alpha=\beta=0$. Therefore, $\psi^*_0(W_3) = (0,0,0)$.
\end{proof}

\subsection*{Acknowledgements}
I am grateful to I.Krichever, B.Khesin, and anonymous referees for important comments and their help in improving this paper.
This work was partially supported by the NSERC research grant.

\end{document}